\documentclass[12pt]{article}

\usepackage{a4wide}
\usepackage{amssymb,amsmath,amsthm}

\usepackage{rotating}
\usepackage{graphicx}
\usepackage{latexsym}
\usepackage{bm}
\usepackage{mathrsfs}

\newcommand{\AGf}{\ensuremath{\mathcal{A}}}
\newcommand{\PGf}{\ensuremath{\mathcal{P}}}
\newcommand{\QGf}{\ensuremath{\mathcal{Q}}}
\newcommand{\PPR}{\ensuremath{\mathcal{P}_{\rm Rot}}}
\newcommand{\PP}[1]{\ensuremath{\mathcal{P}^{(#1)}}}
\newcommand{\PPM}[2]{\ensuremath{\mathcal{P}^{(#1)}_{#2}}}
\newcommand{\StM}[1]{\ensuremath{\mathcal{P}}^{\,\Box\,(#1)}}

\newcommand{\StMM}[2]{\ensuremath{\mathcal{P}}_{#2}^{\,\Box\,(#1)}}
\newcommand{\Pol}{\ensuremath{\mathscr{P}}}
\newcommand{\Sign}[2]{\ensuremath{{ \rm Sign}\left  ( [x^{#2}] #1 \right )}}
\newcommand{\La}[2]{\ensuremath{A^{(#1)}_{#2}}}
\newcommand{\Lb}[2]{\ensuremath{B^{(#1)}_{#2}}}
\newcommand{\Lc}[2]{\ensuremath{C^{(#1)}_{#2}}}

\newcommand{\FM}[1]{\ensuremath{{\cal F}}^{\,\Box\,(#1)}}

\newtheorem{theorem}{Theorem}
\newtheorem{fact}{Fact}
\newtheorem{prop}{Proposition}
\newtheorem{lem}{Lemma}
\newtheorem{definition}{Definition}
\newtheorem{conjecture}{Conjecture}

\title{Area distribution and scaling function \\for punctured polygons}

\author{\sc Christoph Richard\dag, Iwan Jensen\ddag{} and Anthony J. Guttmann\ddag\\
\normalsize
\dag Fakult\"at f\"ur Mathematik, Universit\"at Bielefeld,\\
\normalsize
Postfach 10 01 31, 33501 Bielefeld, Germany\\
\normalsize
\ddag ARC Centre of Excellence for Mathematics and Statistics of Complex Systems,\\ 
\normalsize
Department of Mathematics and Statistics, 
The University of Melbourne,\\ 
\normalsize
Victoria 3010, Australia}

\begin{document}

\maketitle

\centerline{Mathematics Subject Classifications: 05A15, 05A16}
\centerline{PACS numbers: 05.50.+q, 05.70.Jk, 02.10.Ox}

\begin{abstract}
Punctured polygons are polygons with internal holes which are 
also polygons. The external and internal polygons are of the 
same type, and they are mutually as well as self-avoiding. 
Based on an assumption about the limiting area distribution 
for unpunctured polygons, we rigorously analyse the effect of a 
finite number of punctures on the limiting area distribution in 
a uniform ensemble, where punctured polygons with equal perimeter 
have the same probability of occurrence. Our analysis 
leads to conjectures about the scaling behaviour 
of the models.

We also analyse exact enumeration data. For staircase polygons 
with punctures of fixed size, this yields explicit expressions
for the generating functions of the first few area moments.
For staircase polygons with punctures of arbitrary size,
a careful numerical analysis yields very accurate estimates 
for the area moments.
Interestingly, we find that the leading correction term for each area moment 
is proportional to the corresponding area moment with one less puncture. 
We finally analyse corresponding quantities for punctured self-avoiding 
polygons and find agreement with the conjectured formulas to at 
least 3--4 significant digits.
\end{abstract}

\maketitle

\section{Introduction}

The behaviour of planar self-avoiding walks (SAW) and polygons (SAP) is one 
of the classical unsolved problems, not only of algebraic combinatorics, 
but also of chemistry and of physics \cite{MSbook,HughesV1,Vanderzande}.
In the field of algebraic combinatorics, it is a classical enumeration problem.
In chemistry and physics, SAWs and SAPs are used to model a variety of phenomena, including the 
properties of long-chain polymers in dilute solution \cite{deGennes},
the behaviour of ring polymers and vesicles in general \cite{LSF87} and 
benzenoid systems \cite{GCbook,VGJ02} in particular.
Though the qualitative form of the phase diagram \cite{FGW91} is known rigorously,
there is otherwise a paucity of rigorous results. However, there are a few conjectures, 
including the exact values of the critical exponents \cite{N82,N84}, 
and more recently the limit distribution of area and scaling function for SAPs, when 
enumerated by both area and perimeter \cite{RGJ01,C01,R02, RJG03, RJG04}.

Models of planar polygons with punctures arise naturally as cross-sections of
three-dimensional vesicle models. In such cross-sections, there may be holes
within holes, and the number of punctures may be infinite. In this work, we exclude
these possibilities. Whereas our methods can be used to study the former
case, the second situation presents new difficulties, which we have not yet 
overcome\footnote{Since punctured polygons with an unlimited number of punctures have, 
in contrast to polygons without punctures, an (ordinary) perimeter 
generating function with zero radius of convergence \cite{GJO01}, both the phase diagram 
and the detailed asymptotics are clearly going to
be very different from those of polygons without punctures. This is discussed 
further in the conclusion.}.

In this work we consider the effect of a finite number of punctures in polygon
models, in particular we study staircase polygons and self-avoiding polygons on the
square lattice. The perimeter of a punctured polygon \cite{GJWE00, JvRW89} is the 
perimeter of its boundary (both internal and external) while the area of a 
punctured polygon is the area of the enclosed by the external perimeter minus the 
area(s) of any holes\footnote{This has to 
be distinguished from so-called composite polygons \cite{J99}.
The perimeter of a composite polygon is defined as the perimeter of 
the external polygon only, resulting in asymptotic behaviour different
from punctured polygons. Moreover, composite polygons can have more complex 
internal structure than just other polygons.}.
As discussed in section \ref{sec:punc} below, the effect of punctures 
on the critical point and critical exponents of the area 
and perimeter generating function has been the subject of previous studies, but 
the effect of punctures on the critical amplitudes and detailed asymptotics
have not, to our knowledge, been previously considered.

Apart from the intrinsic interest of the problem, we also believe it to be the appropriate 
route to study the detailed asymptotics of {\it polyominoes}, 
since punctured polygons are a subclass of polyominoes.
While we still have some way to go to understand the polyomino phase diagram,
we feel that restricting the problem to this important subclass is the correct route. 

The make-up of the paper is as follows: In the next section we review the known situation 
for the perimeter and area generating functions of punctured polygons and polyominoes. 
In section 3 we review the phase diagram and scaling behaviour of staircase 
polygons and self-avoiding polygons. 
In section 4 we rigorously express the asymptotic behaviour of 
models of punctured polygons in the limit of large perimeter in terms of the 
asymptotic behaviour of the model without punctures, by refining arguments
used in \cite{GJWE00}. This leads, in particular, 
to a characterisation of the limit distribution of the area of punctured polygons. 
This result is then used to conjecture scaling functions of punctured polygons. 
We consider three cases of increasing generality. First, we 
 consider the case of minimal punctures. It is
shown that effects of self-avoidance are asymptotically irrelevant, and that elementary area 
counting arguments yield the leading asymptotic behaviour. We then discuss the 
case of a finite number of punctures of bounded size, and finally the case of a 
finite number of punctures of unbounded size. Results for the latter case
are given for models with a finite critical perimeter generating function
such as staircase polygons and self-avoiding polygons. Whereas the latter two 
cases are technically more involved, the underlying arguments are similar to the case
of minimal punctures. If the critical perimeter generating function of the
polygon model without punctures is finite, then all three cases lead, up to 
normalisation, to the same limit distributions 
and scaling function conjectures.

The next two sections discuss the development and application of 
extensive numerical data to test the results of the previous section. Moreover,
the numerical analysis yields predictions, conjectured to be exact, for the corrections to the
asymptotic behaviour. In particular, section 5 describes the very efficient algorithms 
used to generate the data, while section 6 applies a range of numerical tools to
the analysis of the generating functions for
punctured staircase polygons and then punctured self-avoiding polygons.
Here we wish to emphasise that our work on this problem
involved a close interplay between analytical and numerical work.
Initially, our intention was to check our predictions for scaling functions by
studying amplitude ratios for area moments (given in Table~\ref{tab:ratios}). 
We subsequently discovered numerically the {\em exact} solutions for minimally
punctured staircase polygons. We also obtained very accurate estimates for the 
amplitudes of staircase polygons with one or two punctures of arbitrary size.
From these results we were able to conjecture exact expressions for the amplitudes, 
which in turn spurred us on to further analytical work in order to prove these results.
The final section summarises and discusses our results.

\section{Punctured polygons \label{sec:punc}}

We consider polygons on the square lattice in this article. In particular, 
we study self-avoiding polygons and staircase polygons. 
A self-avoiding polygon on a lattice can be defined as a walk along
the edges of the lattice, which starts and ends at the same lattice point, but 
has no other self-intersections. When counting SAPs, they are generally considered 
distinct up to translations, 
change of starting point, and orientation of the walk, so if there are $p_m$ 
SAPs of length or perimeter $m$ there are $2mp_m$ walks (the factor of two 
arising since the walk can go in two directions). On the square lattice the perimeter
of any polygon is always even so it is natural to count polygons by half-perimeter
instead of perimeter. 
The area of a polygon is the number of lattice cells (times
the area of the unit cell) enclosed by the perimeter of the polygon.
A (square lattice) staircase polygon can be defined as the intersection of two mutually avoiding 
directed walks starting at the same lattice point, moving only to the right or up and 
terminating once the walks join at a vertex. Every staircase polygon is a self-avoiding 
polygon. It is well known that the number $p_m$ of staircase polygons of half-perimeter 
$m$ is given by the $(m-1)^{th}$ Catalan number, $p_m={2m-2 \choose m-1}/m$, with
half-perimeter generating function
\begin{equation}\label{eq:StGf}
\hspace{-2cm}
\PGf(x) =\sum_{m}p_m x^m= \frac{1-2x-\sqrt{1-4x}}{2} \sim \frac{1}{4}-
\frac{1}{2}(1-\mu x)^{2-\alpha} 
\qquad (\mu x\nearrow 1),
\end{equation}
where the connective constant $\mu=4$ and the critical exponent $\alpha=3/2$.
Recall that $f(x)\sim g(x)$ as $x\nearrow x_c$ means that $\lim f(x)/g(x)=1$ 
as $x\to x_c$ from below. In addition, as usual, the rhs is understood as the first
two leading terms in an asymptotic expansion of the lhs about $x=1/\mu$, see 
e.g.~\cite[Sec~1]{BH75}.

\begin{figure}
\begin{center}
\includegraphics[scale=0.7]{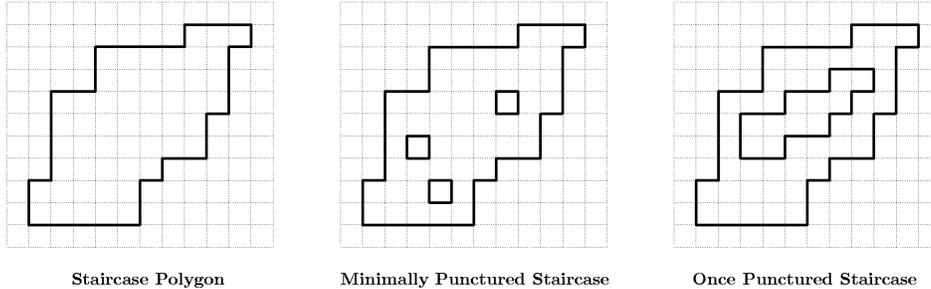}
\end{center}
\caption{\label{fig:poly} 
Examples of the types of staircase polygons we consider in this paper.
}
\end{figure}

Punctured polygons \cite{GJWE00} are polygons with internal holes which are 
also polygons (the polygons are mutually- as well as self-avoiding).  The perimeter 
of a punctured polygon is the sum of the external and internal perimeters
while the area is the area of the external polygon {\em minus} the areas of the
internal polygons. We also consider polygons with {\em minimal} punctures, 
that is, polygons where the punctures are unit cells (or polygons with perimeter 4 
and area 1). Punctured staircase polygons are illustrated in figure~\ref{fig:poly}. 

We briefly review the situation for SAPs with punctures. Analogous results 
can be shown to hold for staircase 
polygons with punctures. Square lattice SAPs with $r$ punctures, counted by area $n$, 
were first studied by Janse van Rensburg and Whittington
\cite{JvRW89}. They proved the existence of an exponential growth constant 
$\kappa^{(r)}$ satisfying $\kappa^{(r)}=\kappa^{(0)}=\kappa$.
Denoting the corresponding number of SAPs by $a_n^{(r)}$ and assuming 
asymptotic behaviour of the form
\begin{displaymath}
a_n^{(r)} \sim A^{(r)}(\kappa^{(r)})^n n^{\beta_r-1}\qquad (n\to\infty),
\end{displaymath}
Janse van Rensburg proved \cite{JvR92} that $\beta_r = \beta_0 + r$. 
These results of course translate to the singular behaviour of the corresponding 
generating functions, defined by $\AGf^{(r)}(q) = \sum_{n > 0} a_n^{(r)}q^n$.

In \cite{GJWE00} Guttmann, Jensen, Wong and Enting studied square lattice 
SAPs with $r$ punctures counted by half-perimeter $m$. They proved the 
existence of an exponential growth constant $\mu^{(r)}$ satisfying
$\mu^{(r)}=\mu^{(0)}=\mu$. If the corresponding number $p_m^{(r)}$
of SAPs is assumed to behave asymptotically as 
\begin{displaymath}
p_m^{(r)}\sim B^{(r)}(\mu^{(r)})^m m^{\alpha_r-3}\qquad (m\to\infty),
\end{displaymath}
they argued, on the basis of a non-rigorous argument, that 
$\alpha_r = \alpha_0 + \frac32 r$. Their results also translate to the associated
half-perimeter generating function $\PP{r}(x) = \sum_{m > 0} p_m^{(r)}x^m$
correspondingly.

Similar results were obtained for polyominoes enumerated by number of cells 
(i.e. area) with a finite number $r$ of punctures \cite{GJWE00}. It has been 
proved that an exponential growth constant $\tau$ exists independently of $r$, 
which satisfies $4.06258\approx \tau >\kappa \approx 3.97087$, where $\kappa$ 
is the growth constant for SAPs enumerated by area. If the number
$a_n^{(r)}$ of polynominoes of area $n$ with $r$ punctures is assumed to 
satisfy asymptotically
\begin{displaymath}
a_n^{(r)} \sim C^{(r)}(\tau^{(r)})^n n^{\gamma_r-1}\qquad (n\to\infty),
\end{displaymath}
it has been shown that $\gamma_r=\gamma_0+r$ and hence that, if the exponents 
$\gamma_r$ exist, they increase by 1 per puncture. It was further conjectured
on the basis of extensive numerical studies \cite{GJWE00}, that the number
$a_n^{(r)}$ satisfies asymptotically
\begin{displaymath}
a_n^{(r)} \sim  \tau^n n^{r-1}\sum_{i \ge 0} C_i^{(r)}/n^{i}
\qquad (n\to\infty).
\end{displaymath}
Notice the conjecture $\gamma_0=0$ and that the correction terms 
go down by a whole power.

For unrestricted polyominoes, that is to say, with no restriction on the number 
of punctures, it was proved by Guttmann, Jensen and Owczarek \cite{GJO01} that 
the  perimeter generating function has zero radius of convergence.
The perimeter is defined to be the perimeter of the boundary plus the
total perimeter of any holes.
If $p_{m}$  denotes the number of polyominoes, distinct up to a translation,
with  half-perimeter $m$, they proved that $p_{m} = m^{m/4 + {\rm o}(m)}$,
meaning that
$$\lim_{m \to \infty}\frac {\log p_{m}}{m \log {m}} = \frac{1}{4}.$$
An attempt to study the quasi-exponential generating function with coefficients
$r_m = p_m/\Gamma(m/4+1)$ was equivocal.
For that reason, studying punctured  self-avoiding polygons was 
considered a controlled route to attempt to
determine the two-variable area-perimeter generating function of polyominoes.

In passing, we note that in \cite{GJ06b} the exact solution of the perimeter generating 
function for staircase polygons with a staircase hole is conjectured, in the form of 
an $8^{th}$ order ODE. It is not obvious how to extract particular asymptotic information, 
notably critical amplitudes from the solution without numerically integrating the ODE. 
In the following, we will obtain such information by combinatorial arguments, 
which refine those of \cite{GJWE00}.

\section{Polygon models and their scaling behaviour}

We review the asymptotic behaviour of self-avoiding polygons
and staircase polygons following mainly \cite{FGW91}. 
For concreteness, consider the fixed perimeter ensemble where, for 
fixed half-perimeter $m$, each polygon of area $n$ has 
a weight proportional to $q^n$, for some positive real number $q$. 
If $0<q<1$, polygons of large area are exponentially suppressed, 
so that typical polygons should be ramified objects. Since such 
polygons would closely resemble branched polymers, the phase $0<q<1$ 
is also referred to as the {\it branched polymer phase}. As $q$ 
approaches unity, typical polygons should fill out more, and 
become less string-like. For $q>1$, polygons of small area are 
exponentially suppressed, so that typical polygons should become 
``fat''. Indeed, they resemble convex polygons \cite{OP99a} and
it has been proved \cite{FGW91} that the mean area of polygons of 
half-perimeter $m$ grows asymptotically proportional to $m^2$.
In the extended phase $q=1$, it is numerically very well established that 
the mean area of polygons of half-perimeter $m$ grows 
asymptotically proportionally to $m^{3/2}$. In the branched polymer 
phase $0<q<1$, the mean area of polygons of half-perimeter $m$ is 
expected to grow asymptotically linearly in $m$, compare also 
\cite[Thm 7.6]{J00} and \cite[Ch~IX.6, Ex.~12]{FS06}.

\begin{figure}[htbp]
  \centering
  \includegraphics[width=5cm]{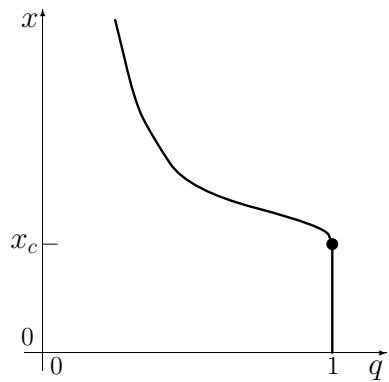}
  \caption{\label{fig:phase} A sketch of the phase diagram of self-avoiding polygons.}
\end{figure}

This change of asymptotic behaviour of typical polygons w.r.t.~$q$ 
is reflected in the singular behaviour of the half-perimeter and area 
generating function
\begin{displaymath}
\PGf(x,q)=\sum_{m,n} p_{m,n} x^m q^n,
\end{displaymath}
where $p_{m,n}$ denotes the number of (self-avoiding) polygons of 
half-perimeter $m$ and area $n$. It has been proved \cite{FGW91} 
that the free energy
\begin{displaymath}
 \kappa (q):= \lim_{m \to \infty} \frac{1}{m} \log \left(\sum_n p_{m,n} q^n\right) 
\end{displaymath}
exists and is finite if $0<q \le 1$. Further, $\kappa(q)$ is log-convex and continuous for these 
values of $q$. It is infinite for $q>1$.
It was proved that for fixed $0<q\le1$, the radius of convergence $x_c(q)$ of $\PGf(x,q)$ 
is given by $x_c(q)= e^{-\kappa(q)}$.
For fixed $q>1$, $\PGf(x,q)$ has zero radius of convergence.
Fisher et al. \cite{FGW91} obtained rigorous upper and
lower bounds on $x_c(q)$.
The expected phase diagram, i.e., the radius of convergence of $\PGf(x,q)$
in the $x-q$ plane, as estimated numerically from extrapolation of SAP 
enumeration data by perimeter and area, 
is sketched qualitatively in figure~\ref{fig:phase}.

For $0<q<1$, the line $x_c(q)$ is, for self-avoiding polygons, expected to be a 
line of logarithmic singularities of the generating function $\PGf(x,q)$. 
For branched polymers in the continuum limit, the existence of the logarithmic 
singularity has recently been proved \cite{BI03}.
The line $q=1$ is, for $0<x<x_c:=x_c(1)$, a line of finite essential singularities 
\cite{FGW91}. For staircase polygons, counted by half-perimeter and area, the 
corresponding phase diagram can be determined exactly, and is qualitatively 
similar to that of self-avoiding polygons. Along the line $x_c(q)$ the 
half-perimeter and area generating function diverges with a simple pole, and 
the line $q=1$ is, for $0<x<x_c$, a line of finite essential 
singularities \cite{Pr94}. 

\medskip 

We will focus on the uniform fixed perimeter ensemble $q=1$
in this article.
Whereas asymptotic area laws in the fixed perimeter ensemble are
expected to be Gaussian for positive $q\ne1$, the behaviour in the uniform 
fixed perimeter ensemble $q=1$ is more interesting. For 
staircase polygons, it can be shown that a limit distribution of 
area exists and is given by the Airy distribution \cite{L84,T91,FL01}. 
For self-avoiding polygons, it is conjectured that an area limit law exists and
is given by the Airy distribution, on the basis of a detailed numerical analysis 
\cite{RGJ01,RJG03,RJG04}. See subsections  \ref{sec:scaling_min} and
\ref{sec:scaling_func}.

If $p_{m,n}$ denotes the number of polygons of half-perimeter $m$ and area $n$,
the existence and the form of a limit distribution can be inferred from the asymptotic
behaviour of the factorial moment coefficients $\sum_{n}(n)_k p_{m,n}$,
where $(a)_k=a\cdot(a-1)\cdot\ldots\cdot(a-k+1)$. The following
result is obtained by standard reasoning \cite{C74}.

\begin{prop}\label{theo:general}
Let for $m,n\in\mathbb N_0$ real numbers $p_{m,n}$ be given.
Assume that the numbers $p_{m,n}$ have the asymptotic form, for $k\in\mathbb N_0$,
\begin{equation}\label{Ecr:ampli}
\sum_n (n)_k p_{m,n}\sim A_k x_c^{-m} 
m^{\gamma_k-1}\qquad (m\to\infty)
\end{equation}
for positive real numbers $A_k$ and $x_c$, where $\gamma_k=(k-\theta)/\phi$, 
with real constants $\theta$ and 
$\phi>0$. Assume that the numbers $M_k:=A_k/A_0$ satisfy
the Carleman condition
\begin{equation}\label{Ecr:growth}
\sum_{k=0}^\infty (M_{2k})^{-1/{2k}}=+\infty.
\end{equation}
Then, for almost all $m$, the random variables $\widetilde X_m$
of area in the uniform fixed perimeter ensemble
\begin{displaymath}
\mathbb P(\widetilde X_m=n)=\frac{p_{m,n}}{\sum_np_{m,n}}
\end{displaymath}
are well defined. We have
\begin{displaymath}
X_m:=\frac{\widetilde X_m}{m^{1/\phi}} \stackrel{d}{\longrightarrow}X \qquad (m\to\infty),
\end{displaymath}
for a uniquely defined random variable $X$ with moments $M_k$,
where the superscript $\stackrel{d}{}$ denotes convergence in distribution.
We also have moment convergence.
\end{prop}

\begin{proof}[Sketch of proof]
A straightforward calculation using Eq.~(\ref{Ecr:ampli}) leads to
\begin{displaymath}
\mathbb E[(\widetilde X_m)_k]
\sim \frac{A_k}{A_0} m^{k/\phi}
\qquad (m\to\infty).
\end{displaymath}
It follows that asymptotically 
the factorial moments are equal to the (ordinary) moments.
Thus, the moments of $X_m$ have the same asymptotic form
\begin{displaymath}
\mathbb E[(X_m)^k]
\sim\frac{A_k}{A_0}=M_k\qquad (m\to\infty).
\end{displaymath}
Due to the growth condition Eq.~(\ref{Ecr:growth}), the sequence 
$(M_k)_{k\in\mathbb N_0}$ defines a unique random variable $X$ with moments $M_k$.
Moment convergence of $(X_m)$ implies convergence in distribution,
see \cite[Thm 4.5.5]{C74} for the line of arguments.
\end{proof}

\medskip

The assumption Eq.~(\ref{Ecr:ampli}) 
translates, on the level of the half-perimeter and area generating function
$\PGf(x,q)$, to a certain asymptotic
behaviour of the so-called factorial moment generating functions
\begin{displaymath}
g_k(x)=\frac{(-1)^k}{k!}\left.\frac{\partial^k}{\partial q^k}
\PGf(x,q)\right|_{q=1}.
\end{displaymath}
It can be shown (compare \cite{F99}) that the asymptotic behaviour Eq.~(\ref{Ecr:ampli}) 
implies for $\gamma_k>0$ the asymptotic equivalence
\begin{equation}\label{eq:relation2}
g_k(x)\sim \frac{f_k}{(x_c-x)^{\gamma_k}}\qquad (x\nearrow x_c),
\end{equation}
where the amplitudes $f_k$ are related to the amplitudes $A_k$%
\footnote{
Note that our definition of the amplitudes $A_k$ differs
from that in \cite{R02} by a factor of $(-1)^kk!$ and from
that in \cite{R05,R06} by a factor of $k!$.} 
in Proposition~\ref{theo:general} via 
\begin{equation}\label{eq:relation}
f_k=\frac{(-1)^k}{k!}A_kx_c^{\gamma_k}\Gamma(\gamma_k).
\end{equation}
If $-1<\gamma_k<0$, the series $g_k(x)$ is convergent as $x\nearrow x_c$, and
the same estimate Eq.~(\ref{eq:relation2}) holds, 
with $g_k(x)$ replaced by $g_k(x)-g_k(x_c)$, where 
$g(x_c):=\lim_{x\nearrow x_c} g(x)$.
In order to deal with these two different cases, we define for 
a power series $g(x)$ with radius of convergence $x_c$, the number
\begin{displaymath}
g^{(c)}=\left\{
\begin{array}{cc}
g(x_c)& {\rm if} \,\,|\lim_{x\nearrow x_c} g(x)|<\infty\\
0 & {\rm otherwise}.
\end{array}
\right.
\end{displaymath}
Adopting the generating function
point of view, the amplitudes $f_k$ determine the numbers $A_k$ and
hence the moments $M_k=A_k/A_0$ of the limit distribution.
The formal series $F(s)=\sum_{k\ge0} f_k s^{-\gamma_k}$ will appear
frequently in the sequel.

\begin{definition}\label{def:fscl}
For the generating function $\PGf(x,q)$ of a class of
self-avoiding polygons, denote its factorial moment 
generating functions by
\begin{displaymath}
g_k(x)=\frac{(-1)^k}{k!}\left.\frac{\partial^k}{\partial q^k}
\PGf(x,q)\right|_{q=1}.
\end{displaymath}
Assume that the factorial moment generating functions satisfy
\begin{equation}\label{form:approx}
g_k(x)-g_k^{(c)} \sim \frac{f_k}{(x_c-x)^{\gamma_k}} \qquad (x\nearrow x_c),
\end{equation}
with real exponents $\gamma_k$. Then, the formal series
\begin{displaymath}
F(s) = \sum_{k\ge0}f_k s^{-\gamma_k}
\end{displaymath}
is called the area amplitude series.
\end{definition}

The area amplitude series is expected to approximate
the half-perimeter and area generating function $\PGf(x,q)$ about 
$(x,q)=(x_c,1)$. This is motivated by the following heuristic 
argument. Assume that $\gamma_k=(k-\theta)/\phi$ with $\phi>0$ and argue
\begin{eqnarray*}
\PGf(x,q)&\approx \sum_{k\ge0} 
\left(g_k^{(c)}+\frac{f_k}{(x_c-x)^{\gamma_k}}\right)(1-q)^k\\
&\approx \left(\sum_{k\ge0}g_k^{(c)}(1-q)^k\right)+(1-q)^{\theta}\left(\sum_{k\ge0}
f_k \left( \frac{x_c-x}{(1-q)^{\phi}}\right)^{-\gamma_k}
\right).
\end{eqnarray*}
In the above calculation, we formally expanded $\PGf(x,q)$
about $q=1$ and then replaced the Taylor coefficients by their
leading singular behaviour about $x=x_c$. 
In the rhs of the above expression, the first sum is by assumption finite,
and the second term contains the area amplitude
series $F(s)$ of combined argument $s=(x_c-x)/(1-q)^\phi$.
This motivates the following definition. A class of self-avoiding
polygons is a subset of self-avoiding polygons. Prominent examples are, 
among others \cite{BM96}, self-avoiding polygons and staircase polygons.

\begin{definition}\label{def:sf}
Let a class of square lattice self-avoiding polygons be given, with
half-perimeter and area generating function $\PGf(x,q)$. 
Let $0<x_c<\infty$ be the radius of convergence 
of the half-perimeter generating function $\PGf(x,1)$. Assume that there 
exist a constant $s_0\in[-\infty,0)$, a function 
${\cal F}:(s_0,\infty)\to\mathbb R$, a real constant $A$ and real
numbers $\theta$ and $\phi>0$, such 
that the generating function $\PGf(x,q)$ satisfies, for real 
$x$ and $q$, where $0<q<1$ and $(x_c-x)/(1-q)^\phi\in(s_0,\infty)$,
the asymptotic equivalence 
\begin{equation}\label{form:scaling}
\PGf(x,q)\sim A+(1-q)^\theta {\cal F}\left(\frac{x_c-x}{(1-q)^{\phi}}\right)
\qquad (x,q)\longrightarrow (x_c,1).
\end{equation}
Then, the function ${\cal F}(s)$  is called a {\it scaling 
function} of combined argument 
$s=(x_c-x)/(1-q)^{\phi}$, and $\theta$ and $\phi$ 
are called {\it critical exponents}.
\end{definition}

\noindent {\bf  Remarks.} {\it i)}
Due to the restriction on the argument of the scaling function, the 
limit $(x,q)\to(x_c,1)$ is approached for values $(x,q)$ satisfying
$x<x_0(q)$ and $q<1$, where $x_0(q)= x_c-s_0(1-q)^\phi$. 
\\
{\it ii)} The above scaling form is also suggested by the theory 
of tricritical scaling, adapted to polygon models \cite{BOP93}. The 
scaling function describes the leading singular behaviour
of $\PGf(x,q)$ about the point $(x_c,1)$ where the two lines of 
qualitatively different singularities meet. \\
{\it iii)} The additional condition $\phi>\theta$ and 
$\theta\notin \mathbb N_0$ ensure that $\gamma_k\in(-1,\infty)
\setminus\{0\}$.
Then, by the above argument, it is plausible that there exists
an asymptotic expansion of the scaling function ${\cal F}(s)$ 
about infinity coinciding with the area amplitude series 
$F(s)$, i.e., ${\cal F}(s)\sim F(s)$ as $s\to\infty$. Recall that
$s$ is considered to be a real parameter.

\medskip

For staircase polygons the existence of a scaling form 
Eq.~(\ref{form:scaling}) has been proved \cite[Thm~5.3]{Pr94},
with scaling function ${\cal F}(s):(s_0,\infty)\to\mathbb R$ 
explicitly given by
\begin{equation}\label{sfstair}
{\cal F}(s) = \frac{1}{16}\frac{\rm d}{{\rm d}s}\log\mbox{Ai} 
\left(  2^{8/3} s \right),
\end{equation}
with exponents $\theta=1/3$ and $\phi=2/3$ and $x_c=1/4$, where
$\mbox{Ai}(x)=\frac{1}{\pi}\int_0^\infty \cos(t^3/3+tx){\rm d}t$ is the
Airy function. The constant $s_0$ is such that $2^{8/3}s_0$ is the
location of the Airy function zero of smallest modulus.  For {\it rooted} 
SAPs with half-perimeter and area generating 
function $\PGf_r(x,q)=x\frac{{\rm d}}{{\rm d}x}\PGf(x,q)$, the conjectured form of 
the scaling function ${\cal F}_r(s):(s_0,\infty)\to\mathbb R$ is \cite{R02}
\begin{displaymath}
{\cal F}_r(s) = \frac{x_c}{2\pi}\frac{\rm d}{{\rm d}s}\log\mbox{Ai} \left( \frac{\pi}{x_c} 
\left( 4A_0\right)^\frac{2}{3} s \right),
\end{displaymath}
with the same exponents as for staircase polygons, $\theta=1/3$ 
and $\phi=2/3$. Here, $x_c=0.14368062927(2)$ is the radius of convergence
of the half-perimeter generating function $\PGf_r(x,1)$ of (rooted) SAPs, and 
$A_0=0.09940174(4)$ is the critical amplitude $\sum_n m p_{m,n}\sim A_0 
x_c^{-m}m^{-3/2}$ of rooted SAPs, which coincides with the critical amplitude
$A_0$ of (unrooted) SAPs. Again, the constant $s_0$ is such that the corresponding 
Airy function argument is the location of the Airy function zero of smallest modulus.
This conjecture was based on the conjecture that both models
have, up to normalisation constants, the same area amplitude series. 
The latter conjecture is supported numerically  
to very high accuracy by an extrapolation of the moment series using exact 
enumeration data \cite{RGJ01,RJG03}. The conjectured form of the scaling 
function ${\cal F}(s):(s_0,\infty)\to\mathbb R$ for SAPs is obtained 
by integration,
\begin{equation}\label{form:unroot}
{\cal F}(s) = -\frac{1}{2\pi} \log\mbox{Ai} 
\left( \frac{\pi}{x_c} \left( 4A_0\right)^\frac{2}{3} s \right)+C(q),
\end{equation}
with exponents $\theta=1$ and $\phi=2/3$. In the above formula, $C(q)$ 
is a $q$ dependent constant of integration, $C(q) = \frac{1}{12 \pi}
(1 - q)\log(1 - q)$, see \cite{RJG04}.
Corresponding results for the triangular and hexagonal lattices can be found in \cite{RGJ01}.

For models of punctured polygons with a finite number of punctures, we have 
qualitatively the same phase diagram as for polygon models without punctures, 
however with different critical exponents $\theta$ depending on the number of 
punctures \cite{JvR92,GJWE00}, and hence we expect different scaling functions. 
We will focus on critical exponents and area limit laws in the uniform
ensemble $q=1$ in the following section. This will lead to conjectures for the
corresponding scaling functions.

\section{Scaling behaviour of punctured polygons \label{sec:scaling}}
We briefly preview the main results of this section.
In subsection~\ref{sec:scaling_min} we study polygons with a finite number
of minimal punctures.  Our result assumes a certain asymptotic form for 
the area moment coefficients for unpunctured polygons. This `assumed' form is 
known to be true for staircase polygons and many other models and universally 
accepted as true for self-avoiding polygons. Given this assumption, we prove 
that the asymptotic behaviour of the area moment coefficients for minimally 
punctured polygons can be expressed in terms of the asymptotic behaviour of 
unpunctured polygons.
In particular we derive expressions for the leading amplitude of the area moments 
for punctured polygons in terms of the amplitudes for unpunctured
polygons. For staircase polygons this leads to exact formulas for the amplitudes.
For self-avoiding polygons the formulas contain certain constants which aren't
known exactly but can be estimated numerically to a very high degree of accuracy.
In subsection~\ref{sec:scaling_bound} we extend the study and proofs to 
polygons with a finite number of punctures of {\em bounded} size and then
in subsection~\ref{sec:scaling_arbitrary} to models with punctures of arbitrary
or unbounded size. Finally in subsection~\ref{sec:scaling_func} we consider
the consequences of our results for the area limit laws of punctured polygons
and we present conjectures for the scaling functions.

\subsection{Polygons with $r$ minimal punctures \label{sec:scaling_min}}

For polygon models with rational perimeter generating functions,
corresponding models with minimal punctures have been studied in 
\cite{RG01}. In particular, 
a method to derive explicit expressions for generating functions
of exactly solvable models with a minimal puncture was 
given \cite[Appendix]{RG01}. It has been applied to Ferrers
diagrams, whose perimeter and area generating function satisfies 
a linear $q$-difference equation, see \cite[Eq.~(54)]{RG01}. 
The method can also be applied to the model of staircase polygons, 
whose half-perimeter and area generating function ${\cal P}(x,q)$
satisfies the quadratic $q$-difference equation
\begin{equation}\label{sp0}
{\cal P}(x,q)=\frac{x^2q}{1-2qx-{\cal P}(qx,q)}.
\end{equation}
Let $\StM{r}(x,q)$ denote the half-perimeter and area generating function
of staircase polygons with $r$ minimal punctures. We have the following 
result for the case $r=1$.
\begin{fact}\label{fact}
The half-perimeter and area generating function of staircase polygons with a 
single minimal puncture $\StM{1}(x,q)$ is given by
\begin{equation}\label{sp1}
\StM{1}(x,q)=\frac{x^4}{(1-2qx-{\cal P}(qx,q))^2} 
\left( {\cal P}(qx,q)-qx \frac{\partial{\cal P} }{\partial x}(qx,q)+
q \frac{\partial {\cal P}}{\partial q}(qx,q) \right),
\end{equation}
where ${\cal P}(x,q)$ satisfies Eq.~(\ref{sp0}).
\qed
\end{fact}
\noindent {\bf Remarks.} 
{\it i)} For a proof of Fact \ref{fact}, proceed along the lines of 
\cite[Appendix]{RG01}. We do not give the details, since we are mainly 
interested in asymptotic results, for which we will give an elementary
combinatorial derivation, valid for arbitrary $r$. See 
Proposition \ref{prop:minpunc} and its subsequent extensions.\\
{\it ii)} For polygons with $r$ punctures, their {\bf $k^{th}$} area moment generating 
functions are defined by $\StM{r}_k(x)= \left(q\frac{\partial}
{\partial q}\right)^k \left.\StM{r}(x,q) \right|_{q=1}$. The above equations can be used 
to obtain explicit expressions for the area moment generating functions $\StM{1}_k(x)$ 
by implicit differentiation. The functions $\StM{1}_k(x)$ also appear in
section \ref{sec:stair_min}.\\
{\it iii)} Assuming that $\StM{1}(x,q)$ has scaling behaviour of the form
\begin{displaymath}
\StM{1}(x,q) \sim (1-q)^{\theta_1}\FM{1}((x_c-x)(1-q)^{-\phi_1})
\end{displaymath}
about $(x,q)=(x_c,1)$, and the necessary analyticity conditions for the
validity of the following calculation, we can 
express the scaling function $\FM{1}(s)$ of staircase polygons with 
a single minimal puncture in terms of the
known scaling function ${\cal F}(s)$ of staircase polygons Eq.~(\ref{sfstair}). 
From Eq.~(\ref{sp1}) we infer that $\theta_1=-2/3$, $\phi_1=2/3$ and
\begin{equation}\label{stsf}
\FM{1}(s)=\frac{1}{24}s{\cal F}'(s)-\frac{1}{48} {\cal F}(s).
\end{equation}

\medskip

In principle, the method of \cite[Appendix]{RG01} can 
be used to analyse the case of several minimal punctures. 
However, the analysis becomes quite cumbersome. On the other 
hand, the previous result suggests simple expressions for 
the scaling functions of models with several punctures in terms 
of that without a puncture. Moreover, we expect such
a phenomenon also to occur for models where an exact 
solution does not exist or is not known. This is 
discussed next. We will asymptotically analyse the area 
moments of a polygon model with punctures and
draw conclusions about their possible scaling behaviour.

For a class of punctured self-avoiding 
polygons, consider their area moment coefficients
\begin{displaymath}
p_m^{\Box(r,k)}:= \sum_n n^kp_{m,n}^{\Box(r)},
\end{displaymath}
where $p_{m,n}^{\Box(r)}$ denotes the number of polygons 
in the class  with $r$ minimal punctures, $r\in\mathbb N_0$, 
of half-perimeter $m$ and area $n$. For simplicity of notation,
we write $p_{m,n}:=p_{m,n}^{\Box(0)}$ and $p_m^{(k)}:=
p_{m}^{\Box(0,k)}$. The area moments in the uniform fixed
perimeter ensemble are expressed in terms of 
the area moment coefficients via
\begin{equation}
\mathbb E[(\widetilde X_m^{\Box(r)})^k]=\frac{\sum_n n^k p_{m,n}^{\Box(r)}}
{\sum_n p_{m,n}^{\Box(r)}} = \frac{p_m^{\Box(r,k)}}{p_m^{\Box(r,0)}}.
\end{equation}

\medskip

\begin{prop}\label{prop:minpunc}
Assume that, for a class of self-avoiding polygons without punctures, 
the area moment coefficients $p_m^{(k)}$ have
the asymptotic form, for $k\in\mathbb N_0$,
\begin{equation}\label{pmas}
p_m^{(k)} \sim A_k x_c^{-m} m^{\gamma_k-1} \qquad (m\to\infty),
\end{equation} 
for numbers $A_k>0$, $x_c>0$ and exponents $\gamma_k=(k-\theta)/\phi$, 
where $\theta$ and $\phi$ are real constants and $0<\phi<1$. 
Then, the area moment coefficient $p_m^{\Box(r,k)}$ of the polygon class
with $r\ge1$ minimal punctures is asymptotically given by, for $k\in\mathbb N_0$, 
\begin{equation}\label{form:minprk}
p_m^{\Box(r,k)} \sim A_k^{(r)} x_c^{-m} m^{\gamma_k^{(r)}-1} \qquad (m\to\infty),
\end{equation}
where $A_k^{(r)}=A_{k+r}x_c^{2r}/r!$ and $\gamma_k^{(r)}=\gamma_{k+r}$. 
\end{prop}

\begin{proof}
We will derive upper and lower bounds on $p_m^{\Box(r,k)}$, which will be 
shown to coincide asymptotically. Let us call two polygons interacting if 
their boundary curves have non-empty intersection. An upper bound is obtained 
by allowing for interaction between all constituents of a punctured polygon. 
Let a polygon $\Pol$ of half-perimeter $m-2r$ and area $n+r$ be given. The 
number of ways of placing $r$ squares inside $\Pol$ is clearly less than
$(n+r)^r/r!$. We thus have
\begin{displaymath}
p_m^{\Box(r,k)}\le {\widetilde p}_m^{(r,k)} := \frac{1}{r!}
\sum_{n\ge 1} n^k (n+r)^rp_{m-2r,n+r} =\frac{1}{r!}
\sum_{n\ge r+1} (n-r)^k n^rp_{m-2r,n}.
\end{displaymath}
By Bernoulli's inequality, we get for ${\widetilde p}_m^{(r,k)}$ the 
bound
\begin{displaymath}
\frac{1}{r!} \sum_{n\ge r+1} \left( n^{k+r}-kr\,n^{k+r-1}\right) 
p_{m-2r,n}\le {\widetilde p}_m^{(r,k)}\le \frac{1}{r!} 
\sum_{n\ge r+1} n^{k+r}p_{m-2r,n}.
\end{displaymath}
For every polygon of perimeter $2s$ and area $t$ we 
have $t\ge s-1$. Thus, for $m$ sufficiently large, we can 
replace the lower bound of summation $r+1$ by zero.
In particular, the latter relation is for $m\ge 3r+2$  
equivalent to
\begin{displaymath}
\frac{1}{r!} \left( p_{m-2r}^{(k+r)}-kr\,p_{m-2r}^{
(k+r-1)} \right)\le {\widetilde p}_m^{(r,k)}\le 
\frac{1}{r!} p_{m-2r}^{(k+r)}.
\end{displaymath}
The assumption Eq.~(\ref{pmas}) on the asymptotic 
behaviour of $p_{m}^{(k)}$ then implies that
\begin{displaymath}
{\widetilde p}_m^{(r,k)}\sim \frac{x_c^{2r}}{r!}p_m^{
(k+r)}\qquad (m\to\infty).
\end{displaymath}
We derive a lower bound by subtracting from the upper bound an 
upper bound on the number of square-square and square-boundary 
interactions. Clearly, square-square interactions are only
present for $r>1$. For a given polygon $\Pol$, the number 
of square-square interactions of $r$ squares is smaller than the number of 
interactions between two squares, where the remaining 
$r-2$ squares may occur at arbitrary positions within the polygon. There
are five possible configurations for an interaction between two
squares, yielding the upper bound $5(n+r)(n+r)^{r-2}$. Thus, the
contribution to ${\widetilde p}_{m}^{(r,k)}$ from square-square
interactions is bounded from above by
\begin{displaymath}
\sum_{n\ge 1} n^k 5(n+r)(n+r)^{r-2}p_{m-2r,n+r}=
5(r-1)!{\widetilde p}_{m}^{(r-1,k)},
\end{displaymath}
which is asymptotically negligible compared to ${\widetilde p}_{m}
^{(r,k)}$. Similarly, the number of configurations arising from
square-boundary interactions is bounded from above by $\sum_{j=1}^r 
4^j(m-2r)^j(n+r)^{r-j}$. This bound is obtained by estimating the 
number of configurations of $j$ squares at the boundary by 
$4^j(m-2r)^j$, the factor 4 arising from edge and vertex interactions,
the factor $(n+r)^{r-j}$ accounting for arbitrary positions of the
remaining $(r-j)$ squares. We thus get an upper bound
\begin{eqnarray*}
\sum_{j=1}^r 4^j (m-2r)^j (r-j)! {\widetilde p}_m^{(r-j,k)}
\sim \sum_{j=1}^r 4^j x_c^{2r-2j}(m-2r)^j p_m^{(k+r-j)}\\  
 \sim 4 x_c^{2r-2}m\, p_m^{(k+r-1)}
\qquad (m\to\infty).
\end{eqnarray*}
By assumption, the latter bound is asymptotically negligible compared to 
${\widetilde p}_{m}^{(r,k)}$. Thus, the lower bound is asymptotically
equal to the upper bound, which yields the assertion of the proposition.
 \end{proof}

\noindent {\bf Remarks.} 
{\it i)} Proposition \ref{prop:minpunc} expresses the asymptotic
behaviour of the area moment coefficients of minimally 
punctured polygons in terms of those of polygons without punctures.
The assumption Eq.~(\ref{pmas}) on the growth of the area moment coefficients 
of the model without punctures is satisfied for the
usual polygon models \cite{BM96}. The asymptotic behaviour of some models 
satisfying $\phi=1$, to which Proposition \ref{prop:minpunc} does not 
apply, has been studied in \cite{RG01}.\\
\noindent {\it ii)}
As discussed in the previous subsection, the amplitudes 
$A_k$ are related to the amplitudes $f_k$ of Eq.~(\ref{form:approx}) by
Eq.~(\ref{eq:relation}), if $\gamma_k\in(-1,\infty)\setminus\{0\}$. 
For staircase polygons, where $\theta=1/3$ 
and $\phi=2/3$, we have explicit expressions for the
amplitudes $A_k$. More generally, it has been shown \cite{R02,R05,R06}
that, for classes of polygon models whose generating function satisfies a 
$q$-functional equation with a square root as the dominant singularity
of their perimeter generating function, we have $f_k=c_kf_1^kf_0^{1-k}$,
where the numbers $c_k$ are, for $k\ge1$, given by
\begin{equation}\label{ck}
\gamma_{k-1}c_{k-1}+\frac{1}{4}\sum_{l=0}^{k}c_{k-l}c_l=0, \qquad c_0=1.
\end{equation}
The critical point $x_c$ as well as $f_0$ and $f_1$ 
are model dependent constants. For staircase polygons we have 
$x_c=1/4$, $f_0=-1$ and $f_1=-1/64$.\\
\noindent {\it iii)} Rooted self-avoiding polygons are conjectured to also have
the exponents $\theta=1/3$ and $\phi=2/3$. In this case the asymptotic form 
Eq.~(\ref{pmas}) and the form of the amplitudes $A_k$, given in 
Eqs.~(\ref{eq:relation}) and (\ref{ck}), has been tested for $k\le10$ 
and shown to hold for  to a high degree of 
numerical accuracy \cite{RJG03}. Here $x_c=0.14368062927(2)$ is the
radius of convergence of the (rooted) SAP half-perimeter generating function,
$f_0=-0.929607(1)$ and $f_1=-x_c/(8\pi)$ are the rooted SAP critical amplitudes
as in Eq.~(\ref{form:approx}). We conjecture that the asymptotic form~(\ref{pmas}) 
holds for rooted SAPs for all values of $k$. Accepting this conjecture to be true, 
Proposition~\ref{prop:minpunc} gives the asymptotic behaviour
for rooted self-avoiding polygons with $r$ minimal punctures.
By definition, unrooted SAPs have the same amplitudes $A_k$.
\\
{\it iv)} The crude combinatorial estimates of interactions in the proof 
of Proposition~\ref{prop:minpunc} cannot be used to obtain corrections to the asymptotic 
behaviour. See also the discussion in the conclusion.

\subsection{Polygons with $r$ punctures of bounded size \label{sec:scaling_bound}}

The arguments in the above proof can be applied to 
obtain results for polygon models with a finite number 
of punctures of bounded size. 
The following theorem generalises Proposition~\ref{prop:minpunc} and serves as  
preparation for the next section, where the case of a finite number 
of punctures of arbitrary size is discussed.
For a class of punctured self-avoiding polygons, consider 
their area moment coefficients
\begin{displaymath}
p_m^{(r,k,s)}:= \sum_n n^kp_{m,n}^{(r,s)},
\end{displaymath}
where $p_{m,n}^{(r,s)}$ denotes the number of polygons in the 
class of half-perimeter $m$ and area $n$ with $r$ punctures,
 $r\in\mathbb N_0$, obeying the condition that the sum of the 
half-perimeter values of the puncturing polygons equals $s$. 
For simplicity of notation, we write $p_{m,n}:=p_{m,n}^{(0,0)}$, 
$p_m^{(k)}:=p_{m}^{(0,k,0)}$ and $p_m:=p_m^{(0)}$.

\begin{theorem}\label{theo:bounded}
Assume that, for a class of self-avoiding polygons without 
punctures, the area moment coefficients $p_m^{(k)}$ have 
the asymptotic form, for $k\in\mathbb N_0$,
\begin{displaymath}
p_m^{(k)} \sim A_k x_c^{-m} m^{\gamma_k-1} \qquad 
(m\to\infty),
\end{displaymath} 
for numbers $A_k>0$,  $x_c>0$ and $\gamma_k=(k-\theta)/\phi$, 
where $\theta$ and $\phi$ are real constants and 
$0<\phi<1$. Denote its half-perimeter generating function 
by ${\cal P}(x)=\left(\sum_{m\ge0} x^m\, p_m\right)$. 
Fix $r\ge1$ und $s\in\mathbb N$ such that $[x^s]({\cal P}(x))^r\ne0$. 
Then, the area moment coefficient $p_m^{(r,k,s)}$ of the polygon class 
with $r\ge1$ punctures whose half-perimeter sum equals is
asymptotically given by, for $k\in\mathbb N_0$,
\begin{equation}\label{eq:Aks}
p_m^{(r,k,s)} \sim A_k^{(r,s)} x_c^{-m}m^{\gamma_k^{(r)}-1}
\qquad (m\to\infty),
\end{equation}
where $\gamma_k^{(r)}=\gamma_{k+r}$ and $A_{k}^{(r,s)}=\frac{A_{k+r}}{r!}
x_c^s[x^s]({\cal P}(x))^r$.
\end{theorem}

\noindent {\bf Remarks.} {\it i)} 
With $s=2$, Theorem~\ref{theo:bounded} reduces 
to Proposition~\ref{prop:minpunc}.
By summation, we also obtain the asymptotic behaviour for models with 
$r$ punctures of total half-perimeter less or equal to $s$. Note that
we have the formal identity 
\begin{displaymath}
\sum_{s=0}^\infty x^s[x^s]({\cal P}(x))^r= ({\cal P}(x))^r.
\end{displaymath}
The above expressions are convergent for $|x|<x_c$. If $\theta>0$, the sum is also
convergent in the limit $x\nearrow x_c$.\\
{\it ii)} The remarks following the proof 
of Proposition~\ref{prop:minpunc} also apply to Theorem~\ref{theo:bounded}.

\begin{proof}
This proof is a direct extension of the proof of Proposition~\ref{prop:minpunc} to the 
case of a finite number of punctures of bounded size. We consider a model of 
punctured polygons where, for fixed $s$, the 
$r$ punctures of half-perimeter $s_i$ and area $t_i$ satisfy 
$s_1+\ldots+s_r= s$. 
We give an asymptotic estimate for $p_m^{(r,k,s)}$. Let a polygon 
$\Pol$ of half-perimeter $m-|\boldsymbol s|$ and of area $n+
|\boldsymbol t|$, where $|\boldsymbol s|=s_1+\ldots+s_r$ and 
$|\boldsymbol t|=t_1+\ldots+t_r$, be given. To obtain an upper 
bound for $p_m^{(r,k,s)}$, ignore all interactions between 
components of a punctured polygon. Recall that two polygons 
interact if their boundary curves have non-empty intersection. 
The number of ways of placing $r$ punctures inside $\Pol$ is 
clearly smaller than 
\begin{displaymath}
(n+|\boldsymbol t|)^r/r!.
\end{displaymath}
This bound is obtained by considering the number of ways of placing 
the lower left corner of each puncture on each square plaquette 
inside the polygon. Note that, unlike in the proof of Proposition 
\ref{prop:minpunc}, this bound also counts configurations where punctures 
protrude from the boundary of $\Pol$. We will compensate for these 
over-counted configurations when deriving a lower bound for 
$p_m^{(r,k,s)}$. We have
\begin{eqnarray}
p_m^{(r,k,s)} &\le {\widetilde p}_m^{(r,k,s)}:= \frac{1}{r!}
\sum_{|\boldsymbol s|= s}\sum_{t_i}\sum_{n\ge 1} n^k (n+|\boldsymbol t|)^rp_{m-
|\boldsymbol s|,n+|\boldsymbol t|}\prod_{i=1}^r p_{s_i,t_i}\nonumber\\
&=\frac{1}{r!} \sum_{|\boldsymbol s|= s}\sum_{t_i}\sum_{n\ge |\boldsymbol t|+1}
(n-|\boldsymbol t|)^k n^r p_{m-|\boldsymbol s|,n}\prod_{i=1}^r \label{form:theo1eq}
p_{s_i,t_i},
\end{eqnarray}
where the first sum is over the variables $s_1,\ldots,s_r$ 
subject to the restriction $|\boldsymbol s|= s$, and the second
sum is over all values of the variables $t_1,\ldots,t_r$.
Note that, for $m$ fixed, all sums are finite. 
Invoking Bernoulli's inequality, we obtain the bound
\begin{eqnarray*}
\frac{1}{r!} \sum_{|\boldsymbol s|= s}\sum_{t_i}&\sum_{n\ge |\boldsymbol t|+1} 
\left( n^{k+r}-k|\boldsymbol t|\, n^{k+r-1}\right) 
p_{m-|\boldsymbol s|,n}\prod_{i=1}^r p_{s_i,t_i}\le 
{\widetilde p}_m^{(r,k,s)}\\
&\le \frac{1}{r!}\sum_{|\boldsymbol s|= s}\sum_{t_i}\sum_{n\ge |\boldsymbol t|+1}
n^{k+r}p_{m-|\boldsymbol s|,n}\prod_{i=1}^r p_{s_i,t_i}.
\end{eqnarray*}
Consider first the asymptotic behaviour of the expression
\begin{displaymath}
{\widetilde a}_{m,s}:= \sum_{|\boldsymbol s|= s}\sum_{t_i} 
\sum_{n\ge |\boldsymbol t|+1} n^{k+r}p_{m-|\boldsymbol s|,n}
\prod_{i=1}^r p_{s_i,t_i}.
\end{displaymath}
If $m\ge |\boldsymbol s|^2+|\boldsymbol s|+2$, then the
lower bound of summation on the index $n$ may be replaced by zero.
This follows from the estimate $t_i\le s_i^2$, being valid
for every self-avoiding polygon of half-perimeter $s_i$ and area $t_i$.
Thus $|\boldsymbol t|\le |\boldsymbol s|^2$, and we argue
that $n\ge m-|\boldsymbol s|-1\ge |\boldsymbol s|^2+1\ge 
|\boldsymbol t|+1$. We thus get for $m$ sufficiently large 
\begin{eqnarray*}
{\widetilde a}_{m,s} = \sum_{|\boldsymbol s|= s} 
p_{m-|\boldsymbol s|}^{(k+r)} \prod_{i=1}^r p_{s_i}
\sim p_m^{(k+r)} \left( \sum_{|\boldsymbol s|= s} 
x_c^{|\boldsymbol s|}\prod_{i=1}^r p_{s_i}\right) \qquad (m\to\infty),
\end{eqnarray*}
where the sum in brackets is finite.
We now analyse the second term in the estimate derived from the
Bernoulli inequality. To this end, define
\begin{displaymath}
{\widetilde b}_{m,s}:= \sum_{|\boldsymbol s|= s}\sum_{t_i} 
\sum_{n\ge |\boldsymbol t|+1} |\boldsymbol t|n^{k+r-1}
p_{m-|\boldsymbol s|,n}\prod_{i=1}^r p_{s_i,t_i}.
\end{displaymath}
Using the estimate $|\boldsymbol t|\le |\boldsymbol s|^2$, we get
\begin{eqnarray*}
{\widetilde b}_{m,s} \le s^2 \left(\sum_{|\boldsymbol s|= s} 
p_{m-|\boldsymbol s|}^{(k+r-1)} \prod_{i=1}^r p_{s_i}\right)
\sim s^2  p_m^{(k+r-1)} \left( \sum_{|\boldsymbol s|= s} 
x_c^{|\boldsymbol s|}\prod_{i=1}^r p_{s_i}\right) \qquad (m\to\infty).
\end{eqnarray*}
Now set  $b_{m,s}:={\widetilde b}_{m,s}/(x_c^{-m}m^{\gamma_{k}^{(r)}-1})$.
The above estimate yields $\lim_{m\to\infty} b_{m,s}=0$, since $0<\phi<1$.

We now derive a lower bound for $p_m^{(r,k,s)}$ by subtracting
from ${\widetilde p}_m^{(r,k,s)}$ an upper bound on the contributions 
arising from puncture-puncture interactions and from puncture-boundary 
interactions. We will show that the lower bound coincides asymptotically
with the upper bound, which then implies the assertion of the theorem
\begin{displaymath}
p_m^{(k,r,s)}\sim \frac{A_{k+r}}{r!} 
\left( \sum_{|\boldsymbol s|= s} x_c^{|\boldsymbol s|}
\prod_{i=1}^r p_{s_i}\right) x_c^{-m}m^{\gamma_{k+r}-1}
\qquad (m\to\infty).
\end{displaymath}
For any polygon $\Pol$, the number of puncture-puncture interactions 
between $r>1$ punctures is smaller than the number of 
puncture-puncture interactions of two punctures with the remaining
$r-2$ punctures occuring at arbitrary positions in the polygon.
We thus get the upper bound
\begin{displaymath}
(t_1+4s_1)t_2(n+|\boldsymbol t|)(n+|\boldsymbol t|)^{r-2}
\le 6 t_1t_2 (n+|\boldsymbol t|)^{r-1},
\end{displaymath}
where we used $t_1\ge s_1-1$. The factor $(t_1+4s_1)t_2$ 
bounds the number of configurations of two interacting punctures,
and the factor $(n+|\boldsymbol t|)^{r-2}$ arises from allowing
arbitrary positions of the remaining $r-2$ punctures. Define
\begin{displaymath}
{\widetilde c}_{m,s}:= \sum_{|\boldsymbol s|= s}\sum_{t_i} 
\sum_{n} t_1t_2 n^k (n+|\boldsymbol t|)^{r-1}
p_{m-|\boldsymbol s|,n+|\boldsymbol t|}\prod_{i=1}^r p_{s_i,t_i}
\le s^4 \left(\sum_{|\boldsymbol s|= s} 
p_{m-|\boldsymbol s|}^{(k+r-1)}\prod_{i=1}^r p_{s_i}\right),
\end{displaymath}
where we used $t_i\le|\boldsymbol t|\le|\boldsymbol s|^2$ for the
last inequality. Setting $c_{m,s}:={\widetilde c}_{m,s}/
(x_c^{-m}m^{\gamma_{k}^{(r)}-1})$, we infer that $\lim_{m\to\infty}
c_{m,s}=0$. We have shown that for fixed $s$ the 
puncture-puncture interactions are asymptotically irrelevant.

We finally estimate the puncture-boundary interactions. This is done
similarly to the above treatment of puncture-puncture interactions.
The number of puncture-boundary interactions is bounded from above by
\begin{displaymath}
\sum_{j=1}^r4^j(m-|\boldsymbol s|)^j s_1\cdot\ldots\cdot s_j \,
(n+|\boldsymbol t|)^{r-j},
\end{displaymath}
where $j$ punctures interact with the boundary, each contributing
a factor $4(m-|\boldsymbol s|)s_i$, and $r-j$ punctures have arbitrary
positions, each contributing a factor $(n+|\boldsymbol t|)$. Note that
the over-counted configurations in ${\widetilde p}_m^{(r,k,s)}$, which 
protrude from the boundary, are compensated for by the above estimate.
Define
\begin{eqnarray*}
{\widetilde d}_{m,s}&:= \sum_{|\boldsymbol s|= s}\sum_{t_i} 
\sum_{n} (m-|\boldsymbol s|)^j n^k (n+|\boldsymbol t|)^{r-j}s_1\cdot
\ldots\cdot s_j\,p_{m-|\boldsymbol s|,n+|\boldsymbol t|}\prod_{i=1}^r 
p_{s_i,t_i}\\
&\le (m-s)^j s^{j} \left(\sum_{|\boldsymbol s|= s} 
p_{m-|\boldsymbol s|}^{(k+r-j)}\prod_{i=1}^r p_{s_i}\right).
\end{eqnarray*}
Defining $d_{m,s}:={\widetilde d}_{m,s}/
(x_c^{-m}m^{\gamma_{k}^{(r)}-1})$, we infer that $\lim_{m\to\infty}
d_{m,s}=0$. We have shown that for fixed $s$ the puncture-boundary 
interactions are asymptotically irrelevant. This completes  the proof. 
\end{proof}

\subsection{Polygons with $r$ punctures of arbitrary size 
\label{sec:scaling_arbitrary}}

For a class of punctured self-avoiding polygons, consider 
for $k\in\mathbb N_0$ their area moment coefficients
\begin{displaymath}
p_m^{(r,k)}:= \sum_n n^kp_{m,n}^{(r)},
\end{displaymath}
where $p_{m,n}^{(r)}:=\sum_{s=0}^\infty p_{m,n}^{(r,s)}<\infty$ denotes 
the number of polygons in the class of half-perimeter $m$ and area $n$ with 
$r$ punctures of arbitrary size, $r\in\mathbb N_0$. For simplicity 
of notation, we write $p_{m,n}=p_{m,n}^{(0)}$ and $p_m^{(k)}=
p_{m}^{(0,k)}$. In the sequel, we will use the area moment generating 
functions $\PGf_k(x)=\sum p_m^{(k)}x^m$ of the model without punctures.

\begin{theorem}\label{theo:arb}
Assume that, for a class of self-avoiding polygons without 
punctures, the area moment coefficients $p_m^{(k)}$ 
have the asymptotic form, for $k\in\mathbb N_0$,
\begin{displaymath}
p_m^{(k)} \sim A_k x_c^{-m} m^{\gamma_k-1} \qquad 
(m\to\infty)
\end{displaymath} 
for numbers $A_k>0$,  $x_c>0$ and $\gamma_k=(k-\theta)/\phi$, 
where $0<\phi<1$. Let $\PGf_k(x)=\sum p_m^{(k)}x^m$ denote the $k^{th}$ area 
moment generating function. 

Then, the area moment coefficient $p_m^{(r,k)}$ of the polygon class with $r\ge1$ 
punctures of arbitrary size is, for $k\in\mathbb N_0$, bounded from above by
\begin{displaymath}
p_m^{(r,k)} \le \frac{[x^m]\PGf_{k+r}(x)(\PGf_0(x))^r}{r!}.
\end{displaymath}
For finite critical perimeter generating functions, characterised by $\theta>0$, 
$p_m^{(r,k)}$ is asymptotically given by, for $k\in\mathbb N_0$, 
\begin{equation}\label{form:t2}
p_m^{(r,k)} \sim \frac{[x^m]\PGf_{k+r}(x)(\PGf_0(x))^r}{r!} \sim
A_k^{(r)}x_c^{-m} m^{\gamma_{k+r}-1}
\qquad (m\to\infty),
\end{equation}
where the amplitudes $A_k^{(r)}$ are given by
\begin{equation}\label{eq:Aka}
A_k^{(r)}=\frac{A_{k+r}(\PGf_0(x_c))^r}{r!},
\end{equation}
where $\PGf_0(x_c):=\lim_{x\nearrow x_c}\PGf_0(x)<\infty$ is the critical amplitude of the
half-perimeter generating function.
\end{theorem}

\noindent {\bf Remarks.} {\it i)}
The asymptotic form Eq.~(\ref{form:t2}) is formally obtained from 
Theorem~\ref{theo:bounded} in the limit of infinite puncture size, 
see Remark {\it i)} after Theorem~\ref{theo:bounded}. This observation 
is also the main ingredient of the following proof, by noting that 
the upper bound has the same asymptotic behaviour.
\\
{\it ii)} For staircase polygons, where 
$\theta=1/3$ and $\phi=2/3$, the assumptions of Theorem 
\ref{theo:arb} are satisfied. For self-avoiding 
polygons, we have the numerically very well established
values $\theta=1$ and $\phi=2/3$, which we believe to describe
the asymptotic behaviour of SAPs.
For models satisfying $\theta<0$, the upper bound generally does
not coincide asymptotically with $p_m^{(r,k)}$. An example
of failure is rectangles with a single puncture.

\begin{proof}
We obtain as in the proof of Theorem~\ref{theo:bounded} an upper bound 
${\widetilde p}_m^{(r,k)}$ for the area moment 
coefficients $p_m^{(r,k)}$. It is given by
\begin{displaymath}
p_m^{(r,k)}\le{\widetilde p}_m^{(r,k)} :=\frac{1}{r!}\sum_{s=0}^m 
\sum_{|{\bf s}|=s}p_{m-|{\bf s}|}^{(k+r)}\prod_{i=1}^r p_{s_i} 
=\frac{1}{r!}[x^m] \PGf_{k+r}(x) (\PGf_0(x))^r.
\end{displaymath}
Assume in the following that $\theta>0$. The asymptotic behaviour 
of the rhs of (\ref{form:t2}) follows by $r$-fold application of 
Lemma 1, which is given in the appendix.
Note that, for $M$ arbitrary, we have by definition 
\begin{displaymath}
p_m^{(r,k)}\ge\sum_{s=0}^M p_m^{(r,k,s)},
\end{displaymath}
where $p_m^{(r,k,s)}$ is the number of $r$-punctured polygons,
whose punctures have total perimeter equal to $s$.
Theorem~\ref{theo:bounded} implies that the above sum is,
for $M$ sufficiently large, asymptotically in $m$, arbitrarily close
to the upper bound ${\widetilde p}_m^{(r,k)}$. See also the 
remark following Theorem~\ref{theo:bounded}.  This 
yields the statement of the theorem.
\end{proof}

\subsection{Limit distribution of area and scaling 
function conjectures \label{sec:scaling_func}}

We first discuss the implications of the previous
results on the asymptotic area law of polygon models
with punctures. By an application of Proposition 
\ref{theo:general}, Theorem \ref{theo:bounded} and 
Theorem \ref{theo:arb} immediately yield the following 
result:

\begin{theorem}\label{thm:limdist}
Assume that, for a class of self-avoiding polygons without 
punctures, the area moment coefficients $p_m^{(k)}$ 
have the asymptotic form, for $k\in\mathbb N_0$,
\begin{displaymath}
p_m^{(k)} \sim A_k x_c^{-m} m^{\gamma_k-1} \qquad 
(m\to\infty)
\end{displaymath} 
for numbers $A_k>0$,  $x_c>0$ and $\gamma_k=(k-\theta)/\phi$, where $0<\phi<1$.
Assume further that the numbers $A_k$ satisfy the Carleman condition
\begin{displaymath}
\sum_{k\ge0} (A_{2k})^{-1/{2k}}=+\infty.
\end{displaymath}
Denote the half-perimeter generating function of the model
by ${\cal P}(x)=\left(\sum_{m\ge0} x^m\, p_m\right)$. 
\begin{itemize}
\item[\it i)] 
Consider for $r\ge1$ the corresponding model with $r$ 
punctures of bounded size, whose half-perimeter sum equals $s\in \mathbb N$, such
that $[x^s]({\cal P}(x))^r\ne0$.
Denote the random variables of area in the uniform fixed 
perimeter ensemble by $\widetilde X_m^{(r,s)}$. Then,
we have convergence in distribution,
\begin{displaymath}
\frac{\widetilde X_m^{(r,s)}}{m^{1/\phi}} \stackrel{d}{\longrightarrow} X^{(r,s)} 
\qquad (m\to\infty),
\end{displaymath}
for a uniquely defined random variable $X^{(r,s)}$ with moments 
\begin{displaymath}
\mathbb E[(X^{(r,s)})^k]=\frac{A_k^{(r,s)}}{A_r},
\end{displaymath}
where the numbers $A_k^{(r,s)}$ are those of Theorem~\ref{theo:bounded}.
We also have moment convergence.
\item[\it ii)]
Let $\widetilde X_m^{(r)}$ denote the
random variable of area in the fixed perimeter ensemble for
the model with $r\ge1$ punctures of unbounded size. If $\theta>0$, then
we have convergence in distribution,
\begin{displaymath}
\frac{\widetilde X_m^{(r)}}{m^{1/\phi}} \stackrel{d}{\longrightarrow} X^{(r)} \qquad (m\to\infty)
\end{displaymath}
for a uniquely defined random variable $X^{(r)}$ with moments 
\begin{displaymath} 
\mathbb E[(X^{(r)})^k]=\frac{A_k^{(r)}}{A_r},
\end{displaymath}
where the numbers $A_k^{(r)}$ are those of Theorem~\ref{theo:arb}.
We also have moment convergence.
\item[\it iii)]
If $\theta>0$, the random variables $X^{(r)}$ and $X^{(r,s)}$ are related by
\begin{displaymath}
x_c^s[x^s]({\cal P}(x))^r\, X^{(r)}=({\cal P}(x_c))^r\, X^{(r,s)},
\end{displaymath}
where ${\cal P}(x)$ is the half-perimeter generating function
of the polygon model without punctures, and where $\mathcal{P}(x_c)=
\lim_{x\nearrow  x_c}\mathcal{P}(x)<\infty$.
\end{itemize}
\qed
\end{theorem}

\noindent {\bf Remarks.} {\it i)}
For a given polygon model satisfying the assumptions of 
Theorem~\ref{thm:limdist}, the area moments satisfy asymptotically
\begin{displaymath}
\frac{\mathbb E[(\widetilde X_m^{(r)})^k]}{k!}\sim D_k^{(r)} m^{k/\phi} 
\qquad (m\to\infty),
\end{displaymath}
for positive numbers $ D_k^{(r)}$.
For classes of polygon models whose generating function satisfies a 
$q$-functional equation with a square root as the dominant singularity
of their perimeter generating function, the amplitude ratios $D_k^{(r)}/\left 
[ D_1^{(r)} \right ]^k$ are universal, i.e., independent of the constants 
$f_0$, $f_1$ and $x_c$, which characterise the underlying model \cite{R02,R05,R06}. This 
follows from Eqs.~(\ref{eq:relation}) and (\ref{ck}) by a straightforward calculation. The 
numbers are listed in Table~\ref{tab:ratios} for small values of $r$. Note that
the same numbers appear for punctures of bounded size. 
\begin{table}
\begin{center}
\begin{tabular}{cccc}
\hline \hline
Amplitude & $r=0$ & $r=1$ & $r=2$ \\
\hline
$D_2/D_1^2$ &
0.530518$\times 10^{-0}$ & 
0.530143$\times 10^{-0}$ &
0.529356$\times 10^{-0}$
\\
$D_3/D_1^3$ & 
0.198944$\times 10^{-0}$ &
0.198369$\times 10^{-0}$ &
0.197361$\times 10^{-0}$
\\
$D_4/D_1^4$ & 
0.592379$\times 10^{-1}$& 
0.588127$\times 10^{-1}$&
0.581533$\times 10^{-1}$
\\
$D_5/D_1^5$ & 
0.149079$\times 10^{-1}$ &  
0.146994$\times 10^{-1}$ &
0.144042$\times 10^{-1}$
\\
$D_6/D_1^6$ &
0.329453$\times 10^{-2}$&
0.321705$\times 10^{-2}$&
0.311511$\times 10^{-2}$ 
\\
$D_7/D_1^7$ &
0.655743$\times 10^{-3}$ & 
0.632288$\times 10^{-3}$ &
0.603260$\times 10^{-3}$
\\
$D_8/D_1^8$ &
0.119654$\times 10^{-3}$& 
0.113600$\times 10^{-3}$&
0.106501$\times 10^{-3}$
\\
$D_9/D_1^9$ & 
0.202754$\times 10^{-4}$ & 
0.189015$\times 10^{-4}$ &
0.173673$\times 10^{-4}$
 \\
$D_{10}/D_1^{10}$ &
0.322150$\times 10^{-5}$ & 
0.294132$\times 10^{-5}$ &
0.264251$\times 10^{-5}$
\\
\hline \hline
\end{tabular}
\end{center}
\caption{\label{tab:ratios} 
Universal amplitude ratios for staircase polygons with $r$ punctures.}
\end{table} 
\\
{\it ii)}
For the above class of models, explicit expressions for the 
asymptotic behaviour of their moment generating functions and 
their probability distributions can be derived from the area 
amplitude series via inverse Laplace transform techniques. 
Since the resulting expressions are quite cumbersome, we do not 
give them here. For $r=0$ the corresponding limit distribution 
of area is the Airy distribution \cite{L84,T91,FL01}. The 
extension to $r\ge1$ is straightforward.
As mentioned above, for $r=0$ the amplitude ratios are found 
to coincide with those of (rooted) self-avoiding polygons to a 
high degree of numerical accuracy \cite{RJG03}. If the conjecture
holds true that they agree {\it exactly}, then the above expressions 
for limit distributions also appear for rooted self-avoiding polygons, 
for all values of $r$. See Section 6 for a detailed numerical 
analysis.

\medskip

We now discuss the relations between the area amplitude series
$F(z)$ of the polygon model without punctures and $F^{(r)}(z)$ of
the polygon model with $r$ punctures. Since all of our models have an
entire moment generating function, the Carleman condition is satisfied,
and Theorem \ref{theo:bounded} and Theorem \ref{theo:arb} yield, by a 
straightforward calculation, the following result.

\begin{theorem}\label{thm:amplfct}
Assume that, for a class of self-avoiding polygons, the polygon model without punctures 
has an area amplitude series, given by
\begin{displaymath}
F(z)=\sum_{k\ge0} f_k z^{-\gamma_k},
\end{displaymath}
where $\gamma_k=(k-\theta)/\phi\in(-1,\infty)\setminus\{0\}$ and $0<\phi<1$, and where
the numbers $f_k\ne0$ are related to the amplitudes $A_k$ in Proposition 
\ref{theo:general} via Eq.~(\ref{eq:relation}). Denote the half-perimeter 
generating function of the model
by ${\cal P}(x)=\sum_{m\ge0} x^m\, p_m$. 
\begin{itemize}
\item[\it i)]
Assume that $r\ge1$ and $s\in\mathbb N$ are given such that
$[x^s]({\cal P}(x))^r\ne0$.
Then, the corresponding model of punctured polygons with $r$ punctures 
of bounded size $s$ has an area amplitude series, given by
\begin{displaymath}
F^{(r,s)}(z)=\sum_{k\ge0} f_k^{(r)} z^{-\gamma_k^{(r)}},
\end{displaymath}
where $\gamma_k^{(r)}=(k-\theta_r)/\phi_r$. We have 
$\theta_r=\theta-r$, $\phi_r=\phi$, and
\begin{equation}\label{form:Fr}
F^{(r,s)}(z)=\frac{(-1)^r}{r!}x_c^s[x^s]({\cal P}(x))^r
\sum_{k\ge r} (k)_r f_k z^{-\gamma_k}.
\end{equation}
\item[\it ii)]
If $\theta>0$, the corresponding model of punctured polygons with $r\ge1$ punctures 
of arbitrary size has an area amplitude series, given by
\begin{displaymath}
F^{(r)}(z)=\frac{(-1)^r}{r!}({\cal P}(x_c))^r
\sum_{k\ge r} (k)_r f_k z^{-\gamma_k},
\end{displaymath}
where $\mathcal{P}(x_c)=\lim_{x\nearrow  x_c}\mathcal{P}(x)<\infty$.
\end{itemize}
\qed
\end{theorem}

\noindent {\bf Remarks.} {\it i)}
Eq.~(\ref{form:Fr}) allows one to derive explicit expressions for the
area amplitude series in terms of $F(z)$. For models where $\theta=1/3$
and $\phi=2/3$ such as staircase polygons,
the area amplitude series $F(z)$ satisfies the Riccati equation
\begin{equation}\label{eq:Ricc}
F(z)^2-4f_1F'(z)-f_0^2z=0.
\end{equation}
This can be used to show that $F^{(r,s)}(z)$ 
(and also $F^{(r)}(z)$) is of the form
\begin{displaymath}
F^{(r,s)}(z)=\sum_{k=0}^{r+1} p_{k,r}(z)F(z)^{k},
\end{displaymath}
where $p_{k,r}(z)$ are polynomials in $z$ of degree not exceeding 
$\lceil 3r/2\rceil$, and $p_{r+1,r}(z)$ is not identically vanishing. 
Simple expressions for the polynomials
$p_{k,r}(z)$ are not apparent.  We note, however, that such
expressions appear as correction-to-scaling functions
of the underlying polygon models without punctures \cite{R02}.\\
{\it ii)}
The model of rooted self-avoiding polygons has been found numerically
to have the same type of area amplitude series as staircase polygons
(with different constants $f_0$ and $f_1$). Similar considerations apply
to the model of unrooted self-avoiding polygons
starting from Eq.~(\ref{form:unroot}).

\medskip

We finally discuss the scaling function conjectures implied
by the results of the previous subsections. For staircase polygons,
the area amplitude
function satisfies the differential equation~(\ref{eq:Ricc}).
This differential equation has a unique solution ${\cal F}(z)$
analytic for $\Re(z)\ge0$, having $F(z)$ as an asymptotic
expansion at infinity, ${\cal F}(z)\sim F(z)$ as $z\to\infty$. 
The function ${\cal F}(z)$ is explicitly given by
\begin{displaymath}
{\cal F}(z)=-4f_1 \frac{{\rm d}}{{\rm d}z}\log\mbox{Ai}
\left(\left(\frac{f_0}{4f_1}\right)^{2/3}z\right),
\end{displaymath}
and this function coincides with the scaling function of
staircase polygons Eq.~(\ref{sfstair}). 

In analogy to the above observation, we conjecture that
the area amplitude series for punctured staircase
polygons determine their scaling functions. Likewise, we conjecture that
the area amplitude series for punctured rooted self-avoiding
polygons determine their scaling functions.

\begin{conjecture}\label{con1}
Let $r\ge1$ and $s\ge2$ be given. For staircase polygons and rooted self-avoiding polygons, 
the area amplitude series $F^{(r,s)}(z)$ and $F^{(r)}(z)$ of
Theorem~\ref{thm:amplfct} uniquely define functions ${\cal F}^{(r,s)}(z)$ 
and ${\cal F}^{(r)}(z)$ analytic for $\Re(z)> z_0$ and non-analytic at $z=z_0$, 
for some negative real number $z_0<0$. We conjecture that the functions 
${\cal F}^{(r,s)}(z):(z_0,\infty)\to\mathbb R$ and 
${\cal F}^{(r)}(z):(z_0,\infty)\to\mathbb R$ are scaling functions
as in  Definition~\ref{def:sf},
\begin{eqnarray*}
{\cal P}^{(r,s)}(x,q) &\sim  (1-q)^{1/3-r}{\cal F}^{(r,s)}
\left(\frac{x_c-x}{(1-q)^{2/3}}\right)\\
{\cal P}^{(r)}(x,q) &\sim  (1-q)^{1/3-r}{\cal F}^{(r)}
\left(\frac{x_c-x}{(1-q)^{2/3}}\right).
\end{eqnarray*}
\end{conjecture}

\medskip

\noindent {\bf Remark.}
The above conjecture has the following implications.\\
{\it i)} Staircase polygons with a single minimal 
puncture  specialise to Eq.~(\ref{stsf}).\\
{\it ii)} Up to constant factors, the scaling 
form of the model with $r$ punctures is obtained as the $r^{th}$ 
derivative w.r.t.~$q$ of the scaling form of the model without 
punctures, as can be proved by induction. As derivatives
can be interpreted combinatorially as marking, this reflects the fact 
underlying the proofs in this section that punctures may be regarded 
as being asymptotically independent, and boundary effects do not 
play a role asymptotically.\\
{\it iii)}
Ignoring questions of analyticity, a (formal) calculation yields 
that the functions ${\cal F}^{(r)}$ (and ${\cal F}^{(r,s)}$) lead, 
for both staircase polygons and (unrooted) self-avoiding polygons, 
to the same critical exponents in the branched polymer phase as 
those conjectured previously \cite{JvR92,GJWE00}. These are obtained from 
the singular behaviour of ${\cal F}^{(r)}$ about the singularity of 
smallest modulus on the negative axis, i.e., at the first zero of 
the Airy-function on the negative axis, see \cite[Sec~1]{R02}. The 
fact that $\PGf^{(r)}(x,q)$ is obtained from $\PGf(x,q)$ by $r$-fold 
differentiation yields the result.

\section{Computer enumerations \label{sec:enum}}

Here we briefly outline which algorithms were used to derive the series expansions
for the area moments of punctured polygons. In most cases (SAPs and punctured
staircase polygons) the algorithms are simple generalisations of previous 
algorithms already described in detail in other papers, referenced below. In
these cases we give brief details of the length of the series and the amount
of CPU time used. Only in the case of staircase polygons with minimal punctures
did we write a new specific algorithm which we shall describe in some detail.

The series for punctured self-avoiding polygons were calculated using a simple
generalisation of the parallel version of the algorithm we used previously 
to enumerate ordinary SAPs \cite{IJ03}. In each case (SAPs with one or two
minimal punctures and SAPs with one or two arbitrary punctures) we calculated the
area moments up to $k=10$ for SAPs to total perimeter 100. Since the smallest such
SAPs have perimeter 16 and 24 this results in series with 42 and 38 non-zero
terms, respectively. The total CPU time required was about 5000 hours for
each of the once punctured SAP problems and up to 3000 hours for the twice
punctured problems. The bulk of these calculations were performed on the old facility
of the Australian Partnership for Advanced Computing (APAC), which
was a Compaq Server Cluster with ES45 nodes with 1GHz Alpha chips (this facility
has since been replaced by an SGI Altix cluster).

In \cite{GJ06b} we used a very efficient algorithm to enumerate once punctured
staircase polygons. The algorithm is easily generalised in order to calculate
area moments which we have done to perimeter 718 ($k=1$), 598 ($k=2$) and 
506 ($k=3$ to 10). It is also quite straight-forward to generalise the algorithm
to count twice punctured staircase polygons and in this case we obtained the
series to perimeter 502 for $k=0$, perimeter 450 for $k=1$ and 2 and to perimeter 
302 for $k=3$ to 10. It is also easy to extract data for staircase polygons with 
punctures of fixed combined perimeter. 
In each case we used around 1000 CPU hours on the APAC Altix cluster which
use 1.6GHz Intel Itanium 2 chips.

\begin{figure}
\begin{center}
\includegraphics[scale=0.6]{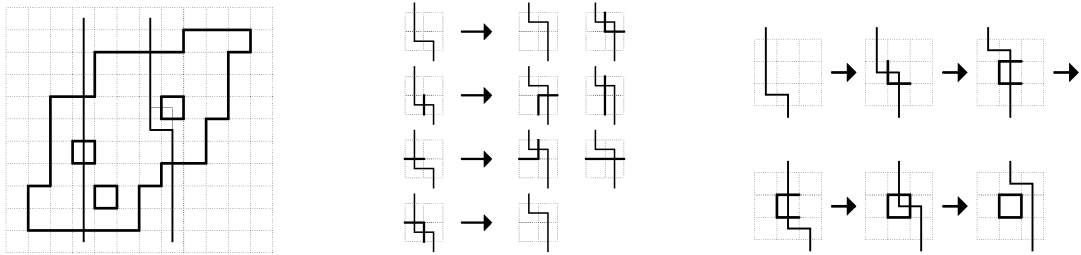}
\end{center}
\caption{\label{fig:TM} 
Illustration of the transfer matrix boundary line and local updating rules.
}
\end{figure}

Finally we describe the algorithm used to enumerate minimally punctured staircase
polygons. The algorithm is based on so-called transfer matrix techniques.
The basic idea is to count the number of polygons by bisecting the lattice with
a boundary line. In the left panel of Fig.~\ref{fig:TM} we show how such a boundary
(the first medium thick line) will intersect the polygon in several places.
The first and last occupied edges intersected by the boundary line are the
directed walks constituting the outer staircase polygon. The other occupied
edges (if any) belong to the minimal punctures. In a calculation to maximal
half-perimeter  $m$ we need only consider intersections with widths up to $w=m/2$.
Any intersection pattern (or signature) can be specified by a string of occupation 
variables, $S = \{\sigma_0,\sigma_1,\ldots \sigma_w \}$, where $\sigma_i=0$, 1 or 2
if edge number $i$ is empty, an occupied outer edge or an occupied edge
part of a minimal puncture, respectively. We could use the same symbol for all
occupied edges but it is convenient to explicitly distinguish between the two cases.
For each signature we keep a generating function which keeps track of the number 
of configurations (to the left of the boundary line), that is, it counts the number 
of possible partially completed polygons with a given signature. In order
to count the total number of punctured polygons we move the boundary line
to the right column by column with each column built up one vertex at a time.
In the left panel of  Fig.~\ref{fig:TM} we have also shown a typical move of the boundary 
line, which starts in the position given by the second medium thick line and where
we add two new edges to the lattice by moving the kink in the boundary line to
the position given by the thin lines. As we move the boundary line to a
new position we calculate the associated generating functions (the updating rules will
be given below). Formally we can view this transformation between signatures as a matrix 
multiplication (hence our use of the nomenclature {\em transfer matrix algorithm}).
However, as can be readily seen, the transfer matrix is extremely sparse and there 
is no reason to list it explicitly (it is given {\em implicitly} by the updating rules).

We start the calculation with the initial signature $\{1,1,0, \ldots ,0\}$, which
corresponds to inserting the two steps of the outer walks in the lower left corner 
(the count of this configuration is 1). As the boundary line is moved it passes over 
a vertex and the updating depends on the states of the edges below and to the left of 
this vertex. After the move we `insert' the edges to the right and above the vertex. 
There  are four possible local configurations of the `incoming' edges as illustrated 
in the middle panels of Fig.~\ref{fig:TM}: Both edges are empty, one of the edges is 
occupied and the other edge is empty or both edges are occupied.

\begin{description}
\item{\bf Both edges empty:} If both incoming edges are empty then both outgoing edges
can be empty. Else we may insert two new steps which {\em must} be part of a minimal puncture
(the outgoing edges are in state `2'). This is {\em only} possible if the vertex is in the
interior of the polygon (there is an edge in state `1' both below and above the vertex).
\item{\bf Left edge empty, bottom edge occupied:} The walk occupying the incoming edge must 
be continued along an outgoing edge. If the occupied edge is part of the external
polygon (in state `1') there are no restrictions. If the occupied edge is part of 
a minimal puncture the walk can only be continued along the edge to the right of the
vertex (otherwise we would not get a minimal puncture).
\item{\bf Left edge occupied, bottom edge empty:} This is similar to the previous
case except that an edge in state `2' must be continued along the edge above the vertex.
\item{\bf Both edges occupied:} If both edges are in state `2' we close the puncture
and the new edges are empty. If the incoming edges are in state `1' we have closed
the outer polygon and then we add the count to the running total for the generating
function.
\end{description}
In the last panel of Fig.~\ref{fig:TM} we show how the updating rules given above
through a sequence of moves of the boundary line gives rise to a minimal puncture.

The perimeter of the completed polygon is given by the position of the boundary
line when the polygon is closed, e.g., if we have taken $k$ steps in the
$x$-direction (completed $k$ columns) and moved the kink $l$ steps in the $y$-direction
then the outer half-perimeter is $k+l$ and the total half-perimeter in $k+l+2r$.
So we need only keep track of the number of punctures $r$. This is done
by associating a truncated polynomial $P_S(x)=\sum p_r x^r$ with each signature, where
the coefficient $p_r$ is the number of partially completed polygons with $r$
punctures of the signature $S$. As a new signature $S'$ is created from $S$
we set $P_{S'}(x)=P_{S'}(x)+x^{\delta}P_S(x)$, where $\delta =1$ if an additional
puncture is inserted (as in the first case described above) or 0 otherwise.
The extension to the calculation of area moments is described in \cite{IJ00}.

\section{Numerical analysis \label{sec:numeric}}

We now turn to our numerical analysis of the series for punctured polygons.
In section~\ref{sec:stair_min} we use our series to determine numerically the 
{\em exact} area moment generating functions for minimally punctured
staircase polygons with up to 5 punctures and $k\leq 10$. The resulting
exact expressions for the leading amplitudes are in complete agreement with
the formula derived in Proposition~\ref{prop:minpunc}. In section \ref{sec:stair_fix} we
extend the study to staircase polygons with one puncture of fixed size and two punctures 
of fixed combined size. Again we find the exact generating functions and 
confirm the formula for the amplitude given in Theorem~\ref{theo:bounded}.
Next, in section~\ref{sec:stair_one}, we analyse area moments for staircase polygons with 
a single puncture of arbitrary size. Guided by results obtained from an analysis of the
conjectured exact ODE \cite{GJ06b} satisfied by the perimeter generating function
we carry out a careful numerical analysis of the area moment coefficients.
This allows to obtain accurate estimates for the leading amplitudes
and we confirm the results of Theorem~\ref{theo:arb} to at least 15 significant digits.
Then, in section~\ref{sec:stair_two} we extend our study to staircase polygons
with two staircase punctures of arbitrary size and we again find good agreement 
with the exact results. Intriguingly, we find in all of the above cases, that
the amplitude of the leading correction term is a constant times the
corresponding amplitude with one less puncture. Finally, in section~\ref{sec:sapana}
we present the results of our analysis of self-avoiding polygons  with one and two punctures 
(minimal as well as arbitrary). In this case the numerical evidence is not 
quite as convincing, but we do find that the numerical estimates agree with the exact formulas
to at least 3--4 significant digits.

\subsection{Staircase polygons with minimal punctures \label{sec:stair_min}}

In \cite{GJWE00} it was found that the half-perimeter generating function of staircase
polygons with a single minimal puncture is: 

\begin{equation}\label{eq:St1p}
\StM{1}(x) = \frac{1-8x+20x^2-16x^3+2x^4}{2(1-4x)}- \frac{1-6x+10x^2-4x^3}{2\sqrt{1-4x}}.
\end{equation}
This result is also derivable from Eq.~(\ref{sp1}), which gives a functional equation 
for the area-perimeter generating function. 
It is thus plausible to expect that the generating function of staircase polygons
with $r$ minimal punctures is of a similar form

\begin{equation}\label{eq:Stmin}
\StM{r}(x) = \frac{A_r(x)+B_r(x)\sqrt{1-4x}}{(1-4x)^{\gamma_r}},
\end{equation}
where $A_r(x)$ and $B_r(x)$ are polynomials and $\gamma_r = (3r-1)/2$. We find this to be correct 
for all the cases we have enumerated that is up to $r=5$:

\begin{eqnarray*}
2A_2(x) &=& {x}^{2}-26{x}^{3}+228{x}^{4}-906{x}^{5}+1709{x}^{6}-1378{x}^{7}+322{x}^{8},\\
2B_2(x) &=& -{x}^{2}+24\,{x}^{3}-182\,{x}^{4}+586\,{x}^{5}-815\,{x}^{6}+404\,{x}^{7}-32\,{x}^{8}.\\
A_3(x) &=& {x}^{2}-22\,{x}^{3}+197\,{x}^{4}-924\,{x}^{5}+2545\,{x}^{6}-5374\,{x}^{7}+13828\,{x}^{8} \\
&&\!-33634\,{x}^{9}\!+46027\,{x}^{10}\!-24746\,{x}^{11}\!+612\,{x}^{12}\!+256\,{x}^{13}\!+256\,{x}^{14},\\
B_3(x) &=& -{x}^{2}+20\,{x}^{3}-159\,{x}^{4}+642\,{x}^{5}-1509\,{x}^{6}+3176\,{x}^{7}-9040\,{x}^{8}\\
&&+19254\,{x}^{9}-18943\,{x}^{10}+4968\,{x}^{11}+768\,{x}^{12}+256\,{x}^{13}.\\
2A_4(x) &=& 2\,{x}^{2}-60\,{x}^{3}+809\,{x}^{4}-6564\,{x}^{5}+36321\,{x}^{6}-146436\,{x}^{7}\\
&&+439283\,{x}^{8}-960070\,{x}^{9}+1485167\,{x}^{10}-1823356\,{x}^{11}\\
&&+2728708\,{x}^{12}-4441406\,{x}^{13}+4054296\,{x}^{14}-932228\,{x}^{15}\\
&&-298318\,{x}^{16}-143360\,{x}^{17}+16384\,{x}^{18}-32768\,{x}^{19},\\
2B_4(x) &=& -2\,{x}^{2}+56\,{x}^{3}-701\,{x}^{4}+5266\,{x}^{5}-26987\,{x}^{6}+100694\,{x}^{7}\\
&&-276415\,{x}^{8}+537888\,{x}^{9}-727683\,{x}^{10}+889018\,{x}^{11}\\
&&-1536634\,{x}^{12}+2199158\,{x}^{13}-1289388\,{x}^{14}-47472\,{x}^{15}\\
&&+26880\,{x}^{16}+50176\,{x}^{17}-6144\,{x}^{18}+8192\,{x}^{19}.\\
2A_5(x) &=& 2\,{x}^{2}-76\,{x}^{3}+1343\,{x}^{4}-14776\,{x}^{5}+114384\,{x}^{6}-666240\,{x}^{7}\\
&&+3036602\,{x}^{8}-11071408\,{x}^{9}+32642310\,{x}^{10}-77911156\,{x}^{11}\\
&&+148630330\,{x}^{12}-220310536\,{x}^{13}+250700412\,{x}^{14}\\
&&-250317844\,{x}^{15}+290657417\,{x}^{16}-309183568\,{x}^{17}\\
&&+150313538\,{x}^{18}+21743832\,{x}^{19}-15222464\,{x}^{20}+449152\,{x}^{21}\\
&&-3828224\,{x}^{22}-2844672\,{x}^{23}+974848\,{x}^{24}-819200\,{x}^{25},\\
2B_5(x) &=& -2\,{x}^{2}+72\,{x}^{3}-1203\,{x}^{4}+12506\,{x}^{5}-91510\,{x}^{6}+504084\,{x}^{7}\\
&&-2171612\,{x}^{8}+7467208\,{x}^{9}-20683474\,{x}^{10}+46059704\,{x}^{11}\\
&&-80841764\,{x}^{12}+107986392\,{x}^{13}-111525400\,{x}^{14}\\
&&+114888220\,{x}^{15}-142562573\,{x}^{16}+122527230\,{x}^{17}\\
&&-24478856\,{x}^{18}-17117496\,{x}^{19}-533632\,{x}^{20}-2988544\,{x}^{21}\\
&&-808960\,{x}^{22}+401408\,{x}^{23}-819200\,{x}^{24}.
\end{eqnarray*}

Likewise the generating functions for
the $k$'th area moment, $\StMM{r}{k}(x)$, is also of the form (\ref{eq:Stmin})
\begin{equation}\label{eq:Stminmom}
\StMM{r}{k}(x) = \frac{A_{r,k}(x)+B_{r,k}(x)\sqrt{1-4x}}{(1-4x)^{\gamma_{r+k}}}.
\end{equation}
We find that the degrees of the polynomials $A_{r,k}(x)$ and $B_{r,k}(x)$ are less than $5r+2k$ 
for $r\leq5$ and $k\leq 10$. In particular we have:

\begin{eqnarray}\label{eq:St1p1m}
\StMM{1}{1}(x) =& \frac{1-14x+72x^2-162x^3+145x^4-34x^5+2x^6}{(1-4x)^{5/2}} \nonumber \\
&-\frac{1-12x+50x^2-82x^3+43x^4-4x^5}{(1-4x)^2}.
\end{eqnarray}

From these solutions we then calculate the exact leading amplitudes 
and indeed we find that
\begin{equation}
\La{r}{k}=A_{k,r}(x_c)/\Gamma (\gamma_{r+k})
=\frac{(-1)^{k+r}(k+r)! x_c^{2r}f_{k+r}}{r! x_c^{\gamma_{k+r}}\Gamma(\gamma_{k+r})}.
\end{equation}
in complete agreement with Eq.~(\ref{eq:relation}).

We have also looked at the amplitudes $\Lb{r}{k}=B_{k,r}(x_c)/\Gamma (\gamma_{r+k}-\frac12)$ 
of the correction terms and find, quite remarkably, that they are given simply in terms
of the leading amplitudes with one less puncture:

\begin{equation}\label{eq:StMB}
\Lb{r}{k}=-\frac18 \La{r-1}{k}.
\end{equation}

\subsection{Staircase polygons with staircase punctures of fixed size \label{sec:stair_fix}}

Next we examine the more general case of staircase polygons with punctures of fixed size. In the
case of one puncture we thus look at staircase polygons with a staircase hole of half-perimeter 
$s$ while in the case of two punctures we look at staircase polygons with two staircase holes 
whose half-perimeters sum to $s$. As in the previous section we expect the generating
functions to be of the form
\begin{equation}\label{eq:Stfix}
\PPM{r,s}{k}(x) = \frac{A_{r,s,k}(x)+B_{r,s,k}(x)\sqrt{1-4x}}{(1-4x)^{\gamma_{r+k}}}.
\end{equation}
We find this to be true with the degree of the polynomials less than $5r+2(k+s)$.

For once punctured polygons we calculated the generating functions for $s \leq 25$ and $k \leq 10$.
In Theorem~\ref{theo:bounded} we proved that the leading amplitude 
$\La{1,s}{k}= A_{1,s,k}(x_c)/\Gamma(\gamma_{k+1})= A_{k+1}x_c^s p_s$,
where $p_s$ is the number of staircase polygons of half-perimeter $s$. This formula is naturally
confirmed by our numerical results. Of more interest is the sub-leading amplitude $\Lb{1,s}{k}$.
We find that the result for minimally punctured polygons generalises to this case and
$\Lb{1,s}{k} = -b_{1,s} A_k$. We also find that the integer sequence 
$d_s=8^{s-1}b_{1,s}=1,5,29,182,\ldots$ is given by the recurrence:
$$8s^2d_s+(s+3)(7s+10)d_{s+1}-(s+3)(s+2)d_{s+2}=0, d_1=0, d_2=1.$$
We note that $d_s$ grows like $8^s$ so that $b_{1,s}$ grows no faster than a polynomial in $s$.

For twice punctured polygons we calculated the generating functions for $s \leq 10$ and $k \leq 10$.
As a consequence of theorem~\ref{theo:bounded} the leading amplitude is given by
$\La{2,s}{k}= A_{2,s,k}(x_c)/\Gamma(\gamma_{k+2})= A_{k+2}x_c^s \sum_{t=2}^{s-2}p_{s-t}p_t$.
And again we find that $\Lb{2,s}{k} = -b_{2,s} A_{k+1}$, though in this case we haven't
found a recurrence for the integer sequence $2\times 8^{s-2}b_{2,s}=1,9,69,510, \ldots$.

\subsection{Staircase polygons with a single staircase puncture \label{sec:stair_one}}

We now turn to the analysis of the  area moments of staircase polygons with a single staircase 
puncture of arbitrary size (1-punctured staircase polygons for short).
In a recent paper \cite{GJ06b} we reported on work which led to an exact Fuchsian linear
differential equation of order 8 apparently satisfied by the half-perimeter generating function, 
$\PP{1}(x) = \sum_{m\geq 0} p_m^{(1)} x^m$, for 1-punctured staircase polygons 
(that is $\PP{1} (x)$ is one of the solutions of the ODE, expanded around the origin). 
Our analysis of the ODE showed that the dominant singular behaviour is
\begin{equation}\label{eq:S1}
\PP{1}(x) \sim \frac{A(x)}{(1-4x)}  + \frac{B(x) + C(x) \log(1-4x)}{\sqrt{1-4x}}+D(x) (1+4x)^{13/2}.
\end{equation}
The functions $A(x)$--$D(x)$ are regular in the disc $|x| \le 1/4$.
In addition there were a pair of singularities on the imaginary axis at $x=\pm i/2$, 
and at the roots of $1+x+7x^2$. Note that the absolute value of these singularities exceeds
$1/4$ and so their contributions to the asymptotic growth of the series coefficients are 
exponentially suppressed.

We expect that the area moment generating functions, $\PPM{1}{k}(x)$, should have a similar 
critical behaviour to that of (\ref{eq:S1}). Indeed our analysis using 
differential approximants \cite{AJG89a} revealed that at $x=x_c=1/4$ there
is a triple root with exponents $-\gamma_{k+1}$ and  $-\gamma_{k+1}\!+\!1/2$ (twice)
which is indicative that the behaviour is
\begin{equation}\label{eq:S1mom}
\PPM{1}{k}(x) \sim \frac{A(x)+[B(x) + C(x) \log(1-4x)]\sqrt{1-4x}}{(1-4x)^{3k/2+1}}.
\end{equation}
However, the behaviour at the singularity $x=x_-=-1/4$ is a little more complicated.
For the first area moment we find that there is a double root with exponents $5$ and
$13/2$, while for the second area moment we find a triple root at $x_-$ with
exponents $7/2$, $5$ and $13/2$. For higher moments the behaviour is consistent 
with a triple root with exponents $(13-3k)/2$, $(10-3k)/2$ and $(7-3k)/2$. That is, the
value of the leading exponent decreases by $3k/2$ and there is in addition a non-analytic 
correction to scaling with exponent $3/2$. The upshot of this analysis is that
the  asymptotic behaviour of the coefficients of $\PPM{1}{k}(x)$ should be given by
\begin{equation}\label{eq:S1co}
 [x^m]\PPM{1}{k}(x)\! \sim\! 4^m\! \left ( \! \sum_{j=0} \!\!
 m^{3k/2-j}\! \left  (\! a_j\!+\!\frac{1}{\sqrt{m}}\left [ b_j \!+\!c_j \log (m) \right ]\! \right ) 
\!\!+\!(-1)^m \! \sum_{j=0}\! d_j m^{(-15+3k-j)/2}\! \right )\!,
\end{equation}
where we have ignored the contributions from singularities in the complex plane with absolute values
exceeding $1/4$. Estimates for the amplitudes were obtained by fitting the coefficients
$p_m^{(1,k)}=[x^m]\PPM{1}{k}(x)$ to the form given above using an increasing number of 
amplitudes. `Experimentally' we find we need about the same total number of terms at 
$x_c$ and $x_-$. So in the fits we used the terms with amplitudes $a_i$, $b_i$ and $c_i$, 
$i=0,\ldots,K$ and $d_i$, $i=0,\ldots,3K$. Going only to $K$ with the $d_i$ amplitudes 
results in much poorer convergence and going beyond $3K$ leads to no improvement. So we 
use the $6K+4$ terms $p_m^{(1,k)}$ with $m$ ranging from $M$ to $M-6K-3$ and solve the
resulting  system of $6K+4$ linear equations. 

\begin{figure}[htbp]
  \centering
  \includegraphics[scale=0.95]{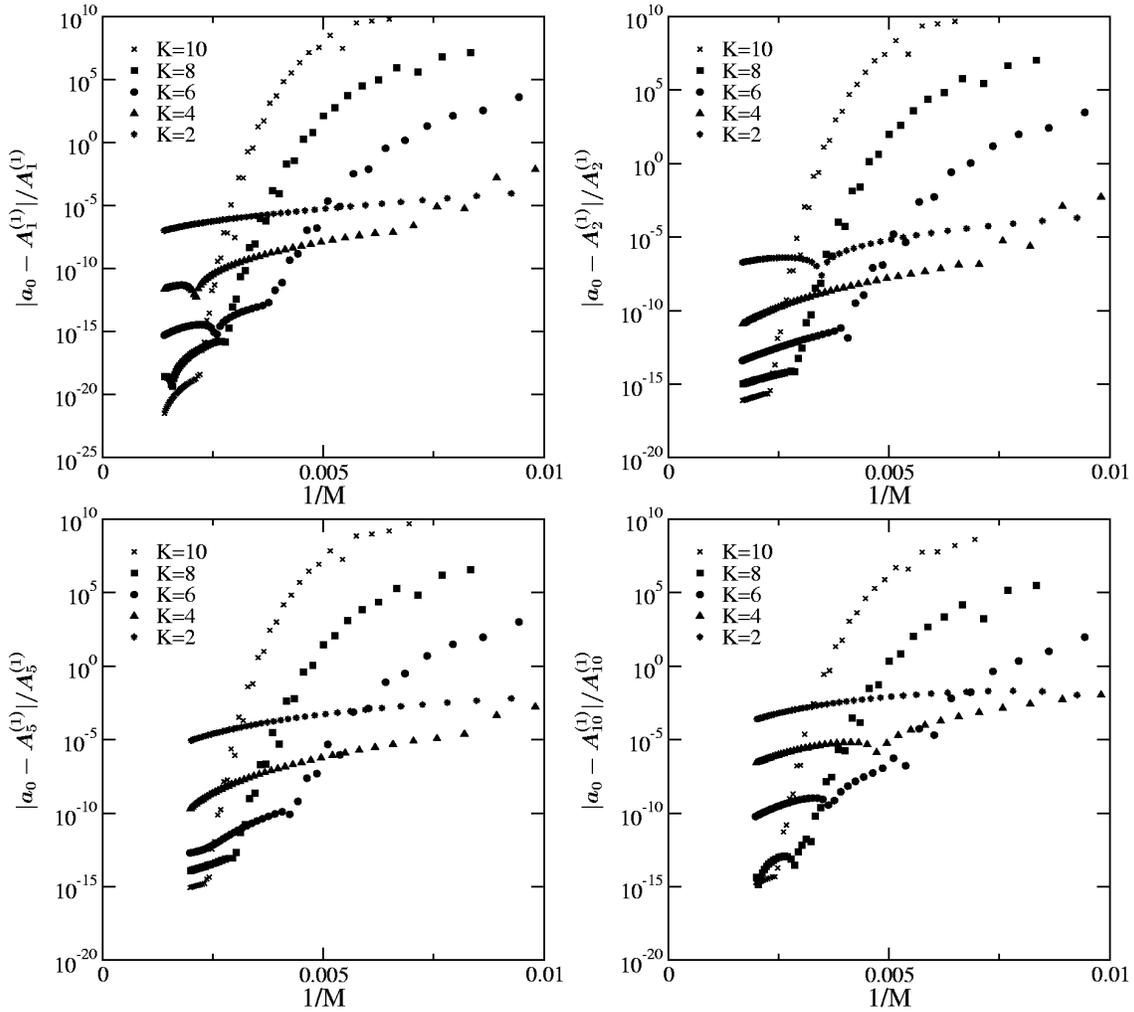}
  \caption{\label{fig:1StA} The relative precision of the estimates for the leading amplitude,
$|a_0 - \La{1}{k}|/\La{1}{k}$, against $1/M$ for the first (top left panel), second
(top right), fifth (bottom left) and tenth (bottom right) area moments.}
\end{figure}

We compare the amplitude estimates to the predictions in Section~\ref{sec:scaling_arbitrary}
and we find that the estimated leading amplitude $a_0$ is given by the exact formula
\begin{equation}\label{eq:StA}
\La{r}{k}=\frac{(-1)^{k+r}(k+r)! x_c^r f_{k+r}}{r! x_c^{\gamma_{k+r}}\Gamma(\gamma_{k+r})},
\end{equation}
which agrees with Eq.~(\ref{eq:Aka}) since the critical half-perimeter generating,
see Eq.~(\ref{eq:StGf}), for staircase polygons is $\PGf(x_c)=1/4=x_c$.
In \cite{GJ06b} we studied the normalised coefficients
$r_m = p_{m+8}^{(1)}/4^m.$ Using the recurrence relations for $p_m^{(1)}$ (derived from the ODE) 
it is easy and fast to generate many more terms $r_m$. We generated the first 100000 terms and saved 
them as floats with 500 digit accuracy. We found (to better than 100 digits) that the leading
amplitude of the normalised series was $\tilde{a}_0=1024$. Going back to the normalisation 
used in this paper we find $a_0=1024/4^8=1/64$ in agreement with formula~(\ref{eq:StA}).
In figure~\ref{fig:1StA} we plot the relative precision of the estimates for the amplitude,
$|a_0 - \La{r}{k}|/\La{r}{k}$ against $1/M$ for some of the area moments.
Recall that for $k=1$ we have a series of 352 non-zero terms, for $k=2$ we have 292 terms
and for $k=3$--10 we have 246 terms. 
In all cases (including for area moments not shown) the relative precision of the estimate 
$a_0$ is better than $10^{-15}$ for $K=10$. For the first area moment, where we have a longer
series, the precision is even more impressive, being better than $10^{-20}$. So in all cases we
can confirm the exact prediction (\ref{eq:StA}) for the amplitude to at least 15 significant
digits.

\begin{figure}[htbp]
  \centering
  \includegraphics[scale=0.95]{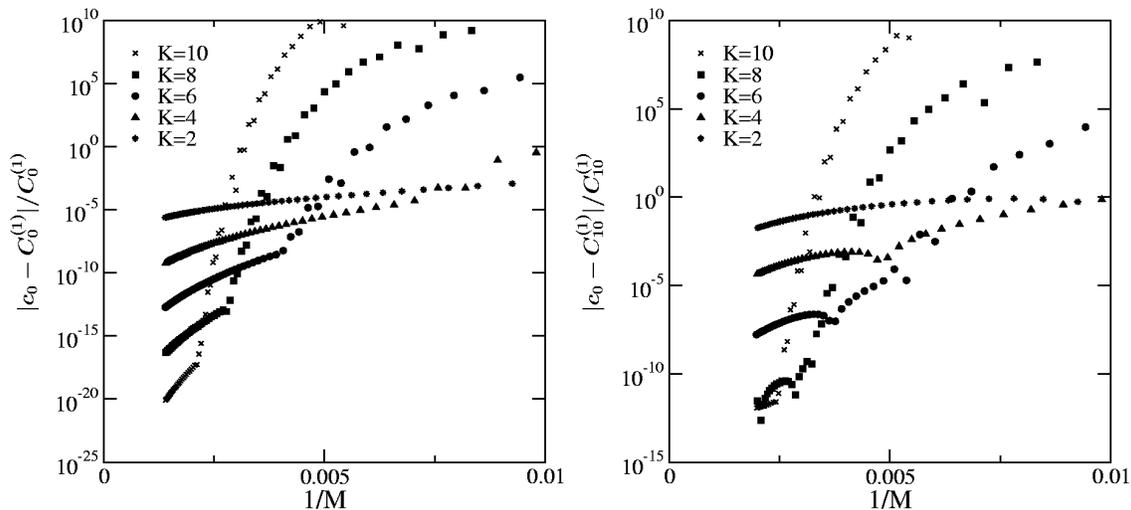}
  \caption{\label{fig:1StC} The relative precision of the estimates for the amplitude,
$|c_0 - \Lc{1}{k}|/\Lc{1}{k}$, against $1/M$ for the zeroth (leftmost panel) 
and tenth (rightmost panel) area moments.}
\end{figure}

In the previous section we saw  that for minimally punctured staircase polygons
the amplitude of the correction term is just a constant times
the leading amplitudes with one less puncture, see Eq.~(\ref{eq:StMB}). It is thus
natural to ask if something similar happens in the more general case.
And indeed we find numerically that the amplitude $\Lc{1}{k}$ of the dominant 
correction in (\ref{eq:S1mom}) (the one proportional to the log-term) giving 
rise to the $c_j$-terms in (\ref{eq:S1co}) are
\begin{equation}\label{eq:S1C}
\Lc{1}{k} = - \frac{3\sqrt{3}}{8\pi} A_k.
\end{equation}
In figure~\ref{fig:1StC} we have plotted the relative error between the estimate
of $c_0$ and the predicted value $\Lc{1}{k}$ for $k=0$ and 10. The estimates for the
other moments are very similar with the accuracy of the agreement diminishing with
higher moments. So in all cases (\ref{eq:S1C}) has been confirmed to better than 10 
digit accuracy. 

For the sub-dominant correction term the amplitude $\Lb{1}{k}$, as approximated
by the term $b_0$ in (\ref{eq:S1co}), is not simply related to $A_k$ and in fact
it even changes sign as $k$ is increased.

\subsection{Staircase polygons with two staircase punctures \label{sec:stair_two}}

In this section we report on our analysis of the series for staircase polygons with
two punctures of arbitrary size. Our first task is to work out the singularity
structure of the perimeter generating function $\PP{2}(x)$ (for which we have
a series with 240 non-zero terms). From the general
considerations of Section~\ref{sec:punc} we expect a singularity at $x=x_c=1/4$ with
exponent $-5/2$, but given the quite complicated singularity structure of  $\PP{1}(x)$,
as detailed in Eq.~(\ref{eq:S1}) and below, we would expect similar complications 
for  $\PP{2}(x)$. We analysed  the series for  $\PP{2}(x)$ using differential
approximants \cite{AJG89a}. This analysis revealed that there is a triple root at $x=x_c=1/4$
and the exponents had values $-2.499(1)$, $-2.070(5)$ and $-1.78(1)$.
So despite having a series of 240 terms it is still very difficult
to pin down the exponents accurately. Given the values quoted above two possible
scenarios present themselves. Either the exponents have the exact values $-5/2$, $-2$ and 
$-2$ or they have the exact values $-5/2$, $-2$ and $-7/4$. The behaviour of
$\PP{1}(x)$ would support the first of these scenarios and below we shall present
evidence from the analysis of the asymptotic form of the coefficients which strongly
supports this behaviour. We also find a double root singularity at $x=x_-=-1/4$
with exponent estimates consistent with the exact values 5 and 11/2. 
In addition there are several conjugate pairs of singularities
in the complex plane. The most important of these are at $x=(-1 \pm i \sqrt{3})/8$,
which has magnitude $1/4$ and thus lies equidistant from the origin to $x_c$. So unlike
the situation for once punctured staircase polygons we cannot ignore the complex singularities.
The exponent estimate at this singularity is consistent with the exact value $33/2$. 
The singularities at $x=(-1 \pm i \sqrt{3})/8$ 
are the roots of the polynomial $1+4x+16x^2$ which would indicate that the generating function 
contains a term $\sim E(x)(1+4x+16x^2)^{33/2}$. Finally, we find some singularities
with magnitude greater than $1/4$. There are singularities at $x=(1 \pm i \sqrt{3})/6$
(which has magnitude $1/3$) with an exponent $33/2$ (we note that these are  the roots
of ($1-3x+9x^2$)) and at $x=\pm i/2$ (magnitude $1/2$) with
an exponent consistent with the value $5$. We also find weak evidence that $\PP{2}(x)$, 
just as $\PP{1}(x),$ has a pair of singularities at the roots of $1+x+7x^2$.

As noted above we need to include the contribution from a conjugate pair of complex 
singularities to the asymptotic form of the coefficients. In general this is not as 
straight-forward a task as dealing with singularities on the real axis. In \cite{IJ06a} 
we examined the generating function of self-avoiding walks on the honeycomb lattice 
which has a pair of singularities on the imaginary axis at $x=\pm i/\tau$, arising 
from a term of the form $H(x)(1+\tau^2 x^2)^{-\eta}$. This typically produces coefficients 
which change sign according to a $++--$ pattern. This can be accommodated by including 
terms of the form

$$
\sim \tau^m m^{\eta-1} \sum_{j\geq 0} (-1)^{\lfloor (m+j)/2 \rfloor} h_j/m^j
$$
in fitting to the asymptotic form of the coefficients.

Note that the analysis in \cite{IJ06a} clearly demonstrated that, as is done above, 
one has to shift the sign-pattern by $j$ when terms proportional to $1/m^j$ are 
considered. Terms arising from other complex conjugate pairs of singularities can 
give rise to much more complicated sign-patterns. In order to handle such cases we 
simply form the Taylor expansion of the simplest possible term arising from the 
singularity and take the sign of the appropriate coefficient. Specifically in 
order to include terms proportional to $1/m^j$, when looking at the coefficients 
$[x^m]\PP{2}(x)$, we take the sign to be the sign of the coefficient of $x^{m+j}$ 
in the Taylor expansion of the function $(1+4x+16x^2)^{33/2}$. We use the notation 
$\Sign{(1\!+\!4x\!+\!16x^2)^{33/2}}{m\!+\!j}$ for this operation.

The singularity structure revealed above leads us to fit the coefficients 
to a form, which is appropriate if the first scenario (exponents at $x=x_c=1/4$ 
are $-5/2$, $-2$ and $-2$) is correct 
\begin{eqnarray}\label{eq:St2Log}
 [x^m]\PP{2}(x)\! =&\! 4^m\! \left (  \sum_{j=0}^K \!
 m^{3/2-j}\! \left  [ a_j\!+\!\frac{1}{\sqrt{m}}\left [ b_j \!+\!c_j \log (m) \right ]\! \right ] 
+(-1)^m \sum_{j=0}^{3K}\! d_j m^{-6-j/2} \right .\nonumber \\
&\left .+ \sum_{j=0}^{3K}\! e_j \Sign{(1\!+\!4x\!+\!16x^2)^{33/2}}{m\!+\!j} m^{-35/2+j}\! \right ),
\end{eqnarray}
where we have ignored the singularities with magnitude exceeding $1/4$. 
We also examine the alternative form appropriate if the second scenario 
(exponents at $x=x_c=1/4$ are $-5/2$, $-2$ and $-7/4$) is correct 
\begin{eqnarray}\label{eq:St2Exp}
 [x^m]\PP{2}(x)\! =&\! 4^m\! \left ( \sum_{j=0}^K \!
 m^{3/2-j}\left [ a_j\!+\!\frac{b_j}{m^{1/2}}\!+\!\frac{c_j}{m^{3/4}}\! \right ] 
 +(-1)^m \sum_{j=0}^{3K}\! d_j m^{-6-j/2}\!\right . \nonumber \\
&\left . + \sum_{j=0}^{3K}\! e_j \Sign{(1\!+\!4x\!+\!16x^2)^{33/2}}{m\!+\!j} m^{-35/2+j}\! \right ). 
\end{eqnarray}

\begin{figure}[htbp]
  \centering
  \includegraphics[scale=0.95]{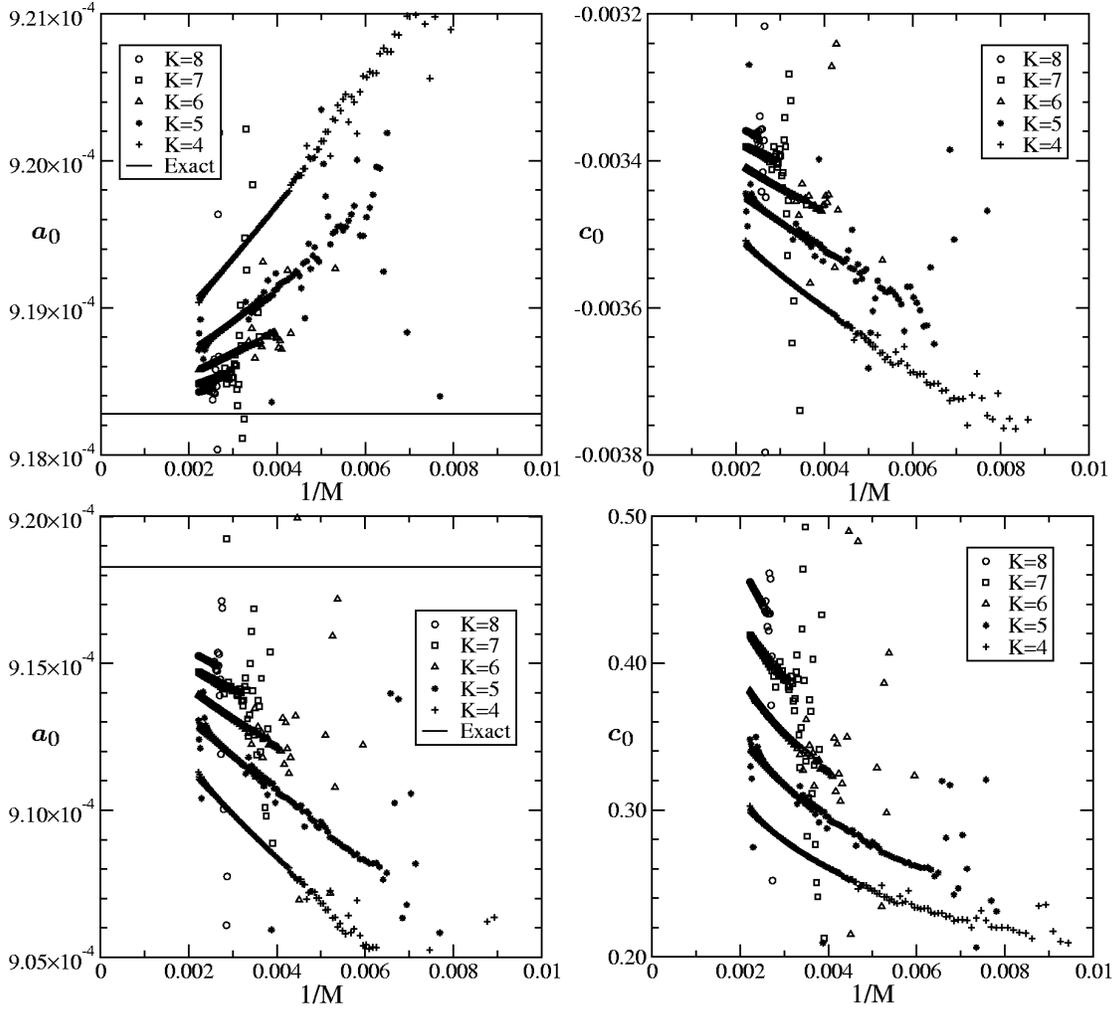}
  \caption{\label{fig:2StA} Estimates of the leading amplitude $a_0$ and the amplitude $c_0$ of the
  sub-dominant term.
  The top panels show the results from fitting to the form~(\protect{\ref{eq:St2Log}})
  while the bottom panels are results from fitting to the form~(\protect{\ref{eq:St2Exp}}).
  The straight line is the exact value of this amplitude $\La{2}{0}$}
\end{figure} 

In figure~\ref{fig:2StA} we have plotted the resulting estimates for the leading amplitude
$a_0$, which we expect are given by the exact value  $\La{2}{0} =\frac{5}{3072\sqrt{\pi}}$, and the 
amplitude $c_0$ of the sub-dominant terms from the two alternative
asymptotic forms. First we focus on the sub-dominant terms. In the top right
panel we show the estimates when fitting to the form~(\ref{eq:St2Log}) where we have a
sub-dominant term proportional to $m \log (m)$ while the bottom right panel shows the
estimates obtained when fitting to the form~(\ref{eq:St2Exp}) where the sub-dominant
term is $m^{3/4}$ (in both cases the dominant term is proportional to $m^{3/2}$ with
a second sub-dominant term proportional to $m$). In the bottom right panel we note that
the estimates for the amplitude $c_0$ of the term $m^{3/4}$ seems to diverge. As $K$
is increased rather than settle down we find that the slope of the curves of the estimates
plotted vs. $1/M$ {\em increases} which is the opposite of what we would expect if we
are fitting to the correct asymptotic form. We take this as firm numerical evidence that 
(\ref{eq:St2Exp}) is incorrect. This could also explain why the corresponding estimates of 
$a_0$ (bottom left panel) don't appear to converge to the predicted exact value. In
contrast the estimates for $c_0$ of the term $m\log (m)$ (top right panel) from the
form~(\ref{eq:St2Log}) do seem to converge and the slopes of the estimates plotted vs. $1/M$ 
decrease with $K$. In this case the estimates for $a_0$ (top left panel) clearly can be
extrapolated to a value consistent with the predicted exact value. Note further
that the top left panel has a resolution along the ordinate which is a factor 5 higher than
in lower right panel so the estimates when fitting to the form~(\ref{eq:St2Log}) are
much more tightly converged. The conclusion is that the numerical evidence clearly
favours the asymptotic form~(\ref{eq:St2Log}) and we believe this to be (if not entirely
correct at least a very good approximation to) the true asymptotic form of the
coefficients of the generating function $\PP{2}(x)$ for twice punctured staircase polygons.

\begin{figure}[htbp]
  \centering
  \includegraphics[scale=0.95]{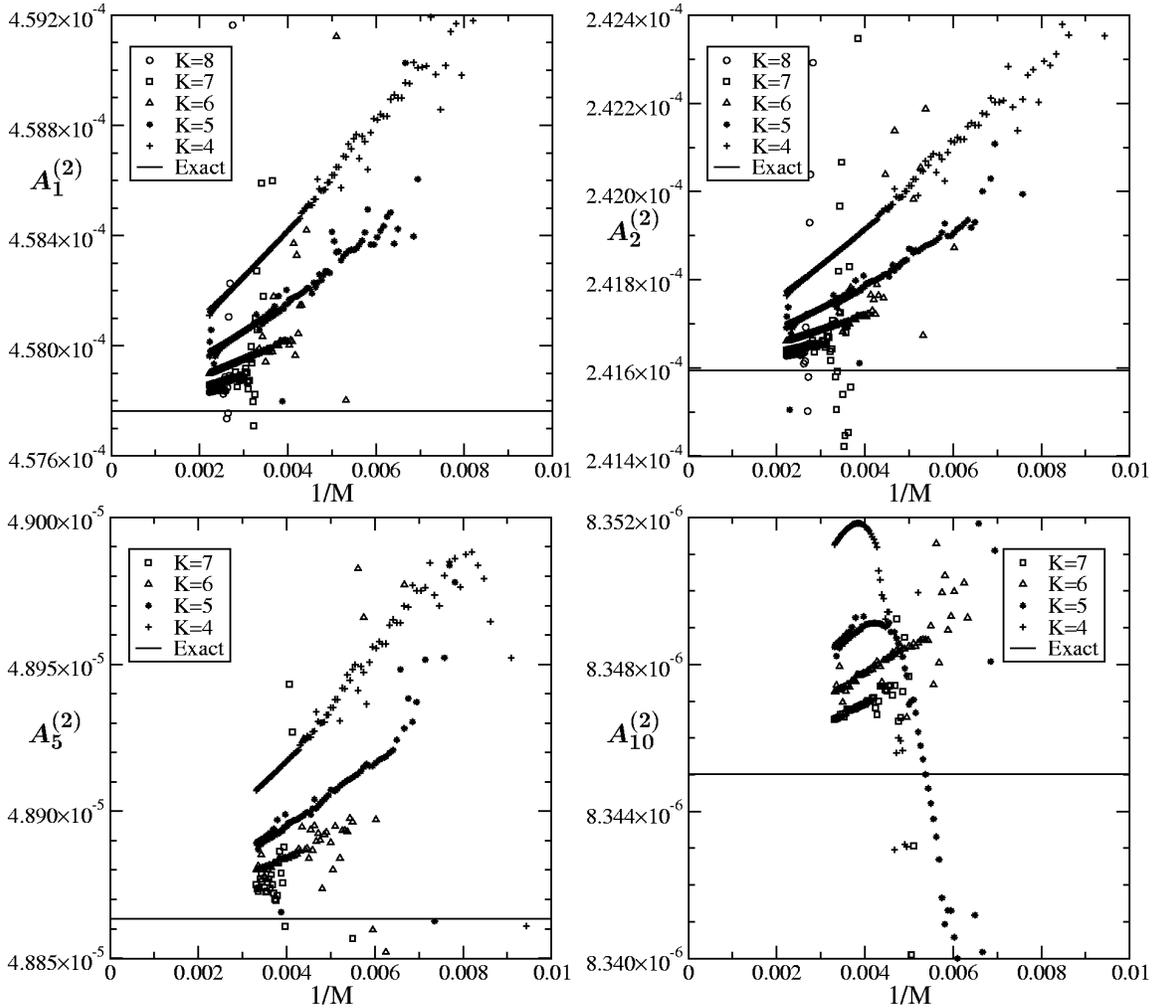}
  \caption{\label{fig:2StM} Estimates of the leading amplitude 
  $\La{2}{k}$ of various area moments
  obtained by fitting to the form~(\protect{\ref{eq:St2Mom}}).
  The straight line is the exact value of this amplitude.}
\end{figure} 

Now that we have settled the question of the correct singularity structure of $\PP{2}(x)$
we turn our attention to the analysis of the area moment generating function $\PPM{2}{k}(x)$.
As for once punctured staircase polygons we find that the leading exponent at all singularities
decreases by $3k/2$. So in order to estimate the leading amplitude $a_0$ we fit to the form
\begin{eqnarray}\label{eq:St2Mom}
 [x^m]\PPM{2}{k}(x)\! =&\! 4^m\! \left (  \sum_{j=0}^K \!
 m^{3(1+k)/2-j}\! \left  [ a_j\!+\!\frac{1}{\sqrt{m}}\left [ b_j \!+\!c_j \log (m) \right ]\! \right ] 
+(-1)^m \sum_{j=0}^{3K}\! d_j m^{(-12+3k+j)/2} \right .\nonumber \\
&\left .+ \sum_{j=0}^{3K}\! e_j \Sign{(1\!+\!4x\!+\!16x^2)^{33/2}}{m\!+\!j} m^{-(35-3k+j)/2}\! \right ),
\end{eqnarray}
where for simplicity our notation suppresses the $k$-dependence of the amplitudes.
Recall that for $k=1$ and 2 we have series with 214 terms and for $k=3$--10 we
have 140 terms.
We compare the amplitude estimates to the predictions of Section~\ref{sec:scaling}.
The leading amplitude $a_0=\La{2}{k}$ is given by the exact 
formula~(\ref{eq:StA}). In figure~\ref{fig:2StM} we show the estimates for the 
leading amplitudes for area moments with $k=1$, 2, 5, and 10.
In all cases the amplitudes estimates appears to converge to the predicted exact
value and agreement is found to at least 3 significant digits.

\begin{figure}[htbp]
  \centering
  \includegraphics[scale=0.95]{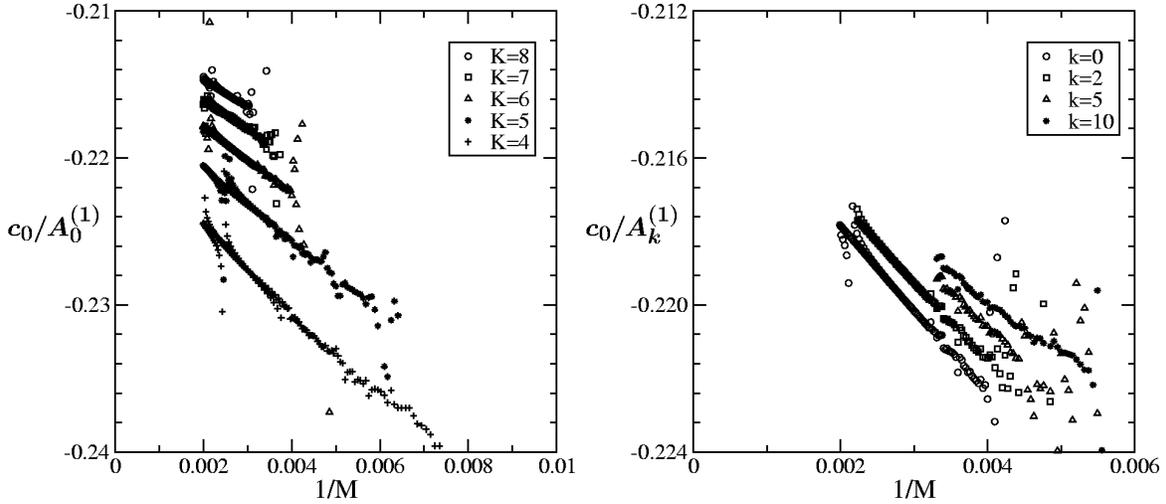}
  \caption{\label{fig:2StC} The ratio $c_0/\La{1}{k}$ for the zeroth area moment 
(leftmost panel) and several different area moments (rightmost panel).}
\end{figure}

Finally we turn our attention to the amplitude $\Lc{2}{k}$ of the dominant correction term.
In the leftmost panel of figure~\ref{fig:2StC} we have plotted the ratio between 
the estimated amplitude $c_0$ and the predicted value of $\La{1}{k}$ for $k=0$ using
several cut-offs. In the rightmost panel we show the same ratio but for
several different moments using the cut-off $K=6$. These plots are consistent
with $\Lc{2}{k} \propto \La{1}{k}$ with the constant of proportionality being
$-0.212(2)$.

\subsection{Punctured self-avoiding polygons \label{sec:sapana}}

\begin{table}
\begin{center}
\small
\begin{tabular}{llll}
\hline\hline
\multicolumn{1}{c}{Degrees} & \multicolumn{3}{c}{Once punctured $k=0$} \\
\hline
$[8,9,9,10]$ &  $ -0.002022  -0.023965i$ & $ -0.002022+  0.023965i$ & $  0.518101$ \\ 
 $[8,9,9,11]$ &  $ -0.003847+  0.034768i$ &  $ -0.003847  -0.034768i$  & $  0.528973$ \\ 
$[8,9,9,12]$ &  $ -0.004553  -0.038104i$ & $ -0.004553+  0.038104i$ & $  0.533388$ \\ 
 $[8,9,10,10]$ &  $ -0.001830  -0.023143i$ & $ -0.001830+  0.023143i$ & $  0.515916$ \\ 
 $[8,9,10,11]$ &  $ -0.044131$ &  $  0.060194$ &  $  0.459558$ \\ 
 $[8,9,11,10]$ &  $ -0.006365  -0.045090i$ & $ -0.006365+  0.045090i$ & $  0.547479$ \\ 
 $[8,10,9,10]$ &  $ -0.007054  -0.047578i$ & $ -0.007054+  0.047578i$ & $  0.552445$ \\ 
 $[8,10,9,11]$ &  $ -0.009398  -0.055135i$ & $ -0.009398+  0.055135i$ & $  0.570859$ \\ 
 $[8,10,10,10]$ &  $ -0.012495  -0.063935i$ & $ -0.012495+  0.063935i$ & $  0.595957$ \\ 
\hline\hline
\end{tabular}
\end{center}
\caption{\label{tab:SAP1pM0} 
Biased estimates of the critical exponents of once punctured SAPs.}
\end{table}

Before proceeding to the estimation of the amplitudes we briefly have a look
at the critical behaviour of the area moment generating functions for
punctured self-avoiding polygons. In \cite{GJWE00} we analysed the behaviour of
$\PPM{r}{k}(x)$  and found  the critical behaviour to be consistent with
$\PPM{r}{k}(x) \sim A(x)+B(x)(x-x_c)^{-\gamma_{k+r}}+ C(x)(x-x_c)^{-\gamma_{k+r}+1/2}$,
where $\gamma_{j} = 3j/2-3/2$.
From previous work \cite{JG99,IJ03} we have very precise estimates for $x_c,$ which
is indistinguishable from the positive root of the polynomial $581x^2+7x-13$, that is,
$x_c=0.1436806292\ldots$. Using this value for $x_c$
we form a $K^{th}$-order {\em biased} differential approximant (DA) to $\PPM{r}{k}(x)$ 
by matching the coefficients in the polynomials $Q_i(x)$ of degree $N_i$
so that (one) of the formal solutions to the homogeneous differential equation
$$
\sum_{i=0}^K (x-x_c)^iQ_i(x)(x\frac{{\rm d}}{{\rm d}x})^i \tilde{P}(x) = 0
$$
agrees with the first $M=\sum_i (N_i+1)$ series coefficients of $\PPM{r}{k}(x)$.
We are thus `forcing' the DAs to have regular singular points at the origin and $x_c$.
The critical exponents $\gamma_j$ ($j=1\ldots K$) are estimated from the indicial equation 
at $x_c$ (note that due to the biasing, $x_c$ is root of order $K$). 
If the true singular behaviour at $x_c$ implies a root of degree less than $K$ we
expect that the `true' exponents will be quite well estimated and show
little scatter while the `surplus' exponents will show a lot of random scatter.
In the following we always use $K=3$ and denote the degrees of the polynomials
$Q_i$ as $[N_3,N_2,N_1,N_0]$.

First we look at the perimeter generating function for once punctured SAPs (for
which we have series with 42 terms).
In this case we have $\gamma_1 =0$ and we thus expect a logarithmic singularity at $x_c$ 
with a square-root correction term as argued in \cite{GJWE00}.
In Table~\ref{tab:SAP1pM0} we list some exponent estimates obtained from the biased
DAs. The exponent estimates are indeed consistent with the exact values $0,0,1/2$,
which confirms the expected behaviour.

\begin{table}
\begin{center}
\small
\begin{tabular}{lllllll}
\hline\hline
\multicolumn{1}{c}{Degrees} & \multicolumn{3}{c}{Once punctured $k=2$} 
                    & \multicolumn{3}{c}{Once punctured $k=5$} \\
\hline
 $[8,9,9,10]$ &  $ -3.000291$ &  $ -2.493921$ &  $ -1.154513$ 
 &  $ -7.500838$ &  $ -6.976201$ &  $ -5.107975$ \\  
$[8,9,9,11]$ &  $ -3.000158$ &  $ -2.496645$ &  $ -1.258606$  
 &  $ -7.503030$ &  $ -6.927432$ &  $ -2.422194$ \\ 
 $[8,9,9,12]$  &  $ -3.012422$ &  $ -2.421400$  &  $ -3.323742$
  &  $ -7.500231$ &  $ -6.992184$ &  $ -5.461825$ \\ 
 $[8,9,10,10]$ &  $ -2.999849$ &  $ -2.502858$ &  $ -1.447159$ 
  & $ -7.517118$ &  $ -6.776728$  &  $ -8.010389$ \\ 
 $[8,9,10,11]$ &  $ -3.000452$ &  $ -2.491160$ &  $ -1.022360$  
  &  $ -7.500700$ &  $ -6.980331$ &  $ -5.251153$ \\ 
 $[8,9,11,10]$ &  $ -3.000143$ &  $ -2.497148$ &  $ -1.307314$  
  &  $ -7.500366$ &  $ -6.988804$ &  $ -5.406296$ \\ 
 $[8,10,9,10]$ &  $ -2.999481$ &  $ -2.511762$ &  $ -1.668527$ 
   &  $ -7.505629$ &  $ -6.879658$ &  $-14.20958$\\ 
 $[8,10,9,11]$ &  $ -2.999906$ &  $ -2.501910$ &  $ -1.441981$  
  &  $ -7.500535$ &  $ -6.984576$ &  $ -5.361394$ \\ 
 $[8,10,10,10]$ &  $ -2.999647$ &  $ -2.507601$ &  $ -1.578352$  
  &  $ -7.500946$ &  $ -6.974159$ &  $ -5.059683$ \\ 
\hline\hline
\end{tabular}
\end{center}
\caption{\label{tab:SAP1pM25} 
Biased estimates of the critical exponents for the 2nd and 5th area moments of once punctured SAPs.}
\end{table}

\begin{table}
\begin{center}
\small
\begin{tabular}{lllrlll}
\hline\hline
\multicolumn{1}{c}{Degrees} & \multicolumn{3}{c}{Twice punctured $k=0$} 
                    & \multicolumn{3}{c}{Twice punctured $k=2$} \\
\hline
 $[7,8,8,9]$ &  $ -1.504250$ &  $ -0.968766$ &  $  6.594627$ 
  &  $ -4.499111$ &  $ -4.006801$ &  $ -2.544930$ \\ 
 $[7,8,8,10]$ &  $ -1.504232$ &  $ -0.968829$ &  $  5.029521$ 
  &  $ -4.495871$ &  $ -4.039199$ &  $ -3.083192$ \\ 
 $[7,8,8,11]$ &  $ -1.504209$ &  $ -0.968910$  &  $-60.90721$
  &  $ -4.494258$ &  $ -4.048789$ &  $ -2.701322$ \\ 
 $[7,8,9,9]$ &  $ -1.504210$ &  $ -0.968909$ &  $ 12.03204$ 
  &  $ -4.494568$ &  $ -4.047486$ &  $ -2.836670$ \\ 
 $[7,8,9,10]$ &  $ -1.504197$ &  $ -0.968951$ &  $ 32.28205$  
  &  $ -4.493872$ &  $ -4.051402$ &  $ -2.672813$ \\ 
 $[7,8,10,9]$ &  $ -1.504205$ &  $ -0.968923$ &  $ 22.27646$  
  &  $ -4.508207$ &  $ -3.915977$ &  $ -1.950665$ \\ 
 $[7,9,8,9]$ &  $ -1.504178$ &  $ -0.969017$ &  $ 18.18298$  
  &  $ -4.495340$ &  $ -4.040712$ &  $ -2.812621$ \\ 
 $[7,9,8,10]$ &  $ -1.504177$ &  $ -0.969020$ &  $ 18.81517$ 
  &  $ -4.495136$ &  $ -4.041237$ &  $ -2.679530$ \\ 
 $[7,9,9,9]$ &  $ -1.504176$ &  $ -0.969022$ &  $ 18.58757$ 
  &  $ -4.498396$ &  $ -4.012876$ &  $ -2.698591$ \\ 
\hline\hline
\end{tabular}
\end{center}
\caption{\label{tab:SAP2pM02} 
Biased estimates of the critical exponents for the 0th and 2nd area moments of twice punctured SAPs.}
\end{table}
 
In Table~\ref{tab:SAP1pM25} we list some exponent estimates for the 2nd and 5th 
area moments of once punctured SAPs. The exponents support the expectations for the
leading and sub-dominant exponent. The third exponent is of no significance--the nature
of the differential approximant forces a third exponent to have some value, but its lack
of convergence suggests it is not, in fact present.  In Table~\ref{tab:SAP2pM02} we list exponent
estimates for the 0th and 2nd area moments of twice punctured SAPs (for
which we have series with 38 terms). Similar comments
apply to the three columns of exponent estimates as were made for the 2nd area moments.
In all cases we get a clear confirmation of the critical behaviour observed in \cite{GJWE00}.

Next we turn our attention to estimates for the critical amplitudes.
Proposition~\ref{prop:minpunc} and Theorem~\ref{theo:arb} tells us that the critical amplitude
of the $k$th area moment of self-avoiding polygons with $r$ (minimal or arbitrary)
punctures is proportional to the critical amplitude of the  $(k+r)$th area moment 
of unpunctured SAPs (the theorems also give the constants of proportionality). 
In order to test these predictions numerically we analyse in this section data
for SAPs with one and two punctures. In all cases we look at the ratio 
$r_m=p_m^{(r,k)}/p_m^{(k+r)}$ which should approach the relevant constant of 
proportionality. Given the critical behaviour outlined above we expect that these 
amplitude ratios can be approximated quite well by the asymptotic form
$$r_m = \sum_{j=0} a_j m^{j/2}$$
So as in the analysis of punctured staircase polygons we obtain estimates for the 
leading amplitude $a_0$, by fitting to the above form truncated at some level $K$, 
and we then plot these estimates against $1/M$.

\subsubsection{Minimal punctures}

According to Proposition~\ref{prop:minpunc} the amplitude of the $k$-th area moment of 
SAPs with $r$ minimal punctures is 

\begin{equation}
A_k^{(r)}=A_{k+r}x_c^{2r}/r!.
\end{equation}

\begin{figure}[htbp]
  \centering
  \includegraphics[scale=1.0]{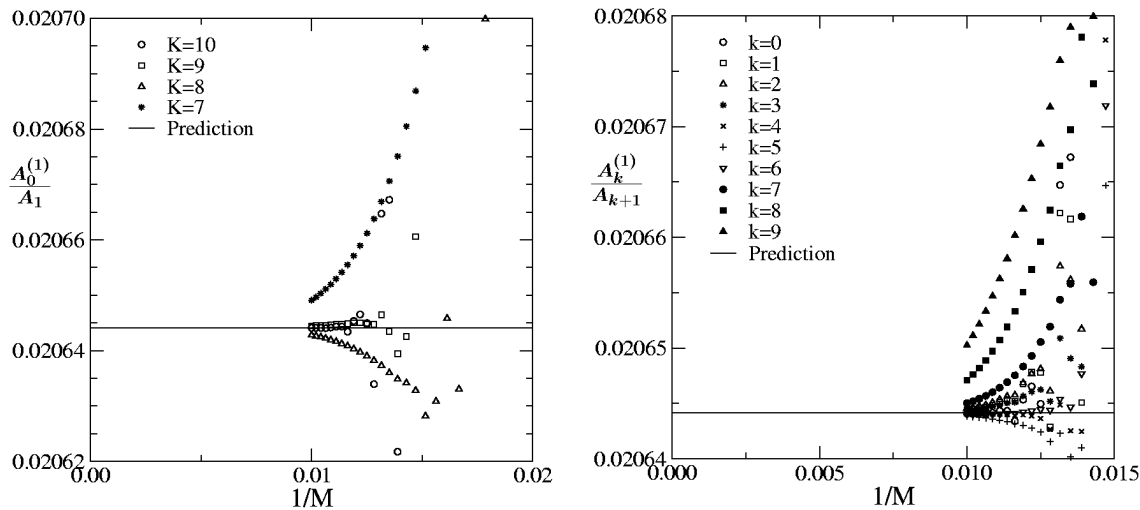}
  \caption{\label{fig:SapM1A} Estimates of the amplitude ratio $A_k^{(1)}/A_{k+1}$ 
   for self-avoiding polygons with one minimal puncture.}
\end{figure}

We first analyse the area moments of SAPs with a single minimal puncture. 
In the left panel of figure~\ref{fig:SapM1A} we have plotted estimates for amplitude ratio
$A_0^{(1)}/A_1$ with $7 \leq K \leq 10$. The prediction for this ratio is $x_c^2$, which
is plotted as a straight line. The estimates obtained with $K=9$ and $10$ are indistinguishable
from the predicted value at the resolution of the plot. Note that the `curvature' of
the plotted values decreases as $K$ increases. We take this to be a very clear indication
that the ratio $r_m$ is very well approximated by the assumed asymptotic form.
In the right panel we plot the estimates for the amplitude ratio $A_k^{(1)}/A_{k+1}$ 
for $0 \leq k \leq 9$ using the cut-off value $K=10$. For small values of $k$ these plots
give firm numerical evidence for the correctness of  Proposition~\ref{prop:minpunc}. For higher
values of $k$ (8 and 9 in particular) the evidence is not quite as firm though nothing
in the plot would suggest a discrepancy with the predicted value. We again emphasise
that in this case we have relatively short series of only 42 terms.

\begin{figure}[htbp]
  \centering
  \includegraphics[scale=1.0]{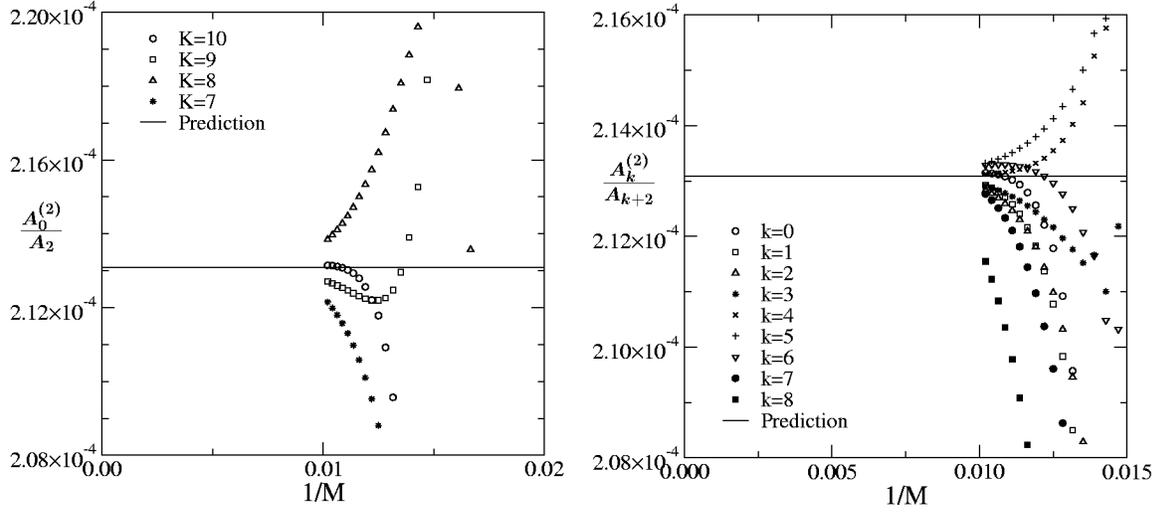}
  \caption{\label{fig:SapM2A} Estimates of the amplitude ratio $A_k^{(2)}/A_{k+2}$ 
  for self-avoiding polygons with two minimal punctures.}
\end{figure}

Next we analyse the area moments of SAPs with two minimal punctures. 
The left panel of figure~\ref{fig:SapM2A} shows estimates for amplitude ratio
$A_0^{(2)}/A_2$ with $7 \leq K \leq 10$. The prediction for the ratio, $x_c^4/2$,
is plotted as a straight line.  Again we find that the estimates for $K=10$ are
indistinguishable from the predicted value, though in this case the relative 
resolution is coarser than in the previous plot. 
In the right panel we plot estimates of the amplitude ratio $A_k^{(2)}/A_{k+2}$ 
for $0 \leq k \leq 8$ using the cut-off value $K=10$. Again we find that our
numerical analysis confirms the prediction to a high degree of confidence.
Recall that in this case we have series with only 38 terms.

\subsubsection{Arbitrary punctures}

According to Theorem~\ref{theo:arb} the amplitude of the $k$-th area moment of 
SAPs with $r$ arbitrary punctures is 

\begin{equation}
A_k^{(r)}=A_{k+r}({\cal P}(x_c))^r/r!, 
\end{equation}
where $  {\cal P}(x_c)$ is the critical amplitude of the half-perimeter generating function.

The first step in our analysis is to obtain an accurate estimate of ${\cal P}(x_c)$.
In \cite{JG99} we obtained the rather imprecise estimate ${\cal P}(x_c)\approx 0.036$ 
by evaluating Pad\'e approximants to the generating function. Here we shall
estimate ${\cal P}(x_c)$ directly from the perimeter data. We first tried the form
\begin{displaymath}
\sum_{m=0}^M p_{m}x_c^m \sim  {\cal P}(x_c)+a_1/M^{1/2}+a_2/M+\cdots
\end{displaymath}
Using the first ten terms in this asymptotic expansion we found (with $M=55$) 
that ${\cal P}(x_c)\approx 0.0362642$, but $a_1\approx 4.28\times10^{-9}$ and
$a_2\approx -1.42\time 10^{-7}$, while $a_3 \approx -0.066$. We are therefore
quite confident that $a_1=a_2=0$. Upon further analysis we found convincing evidence
that the correct asymptotic form in fact is 

\begin{displaymath}
\sum_{m=0}^M p_{m}x_c^m \sim  {\cal P}(x_c)+b_1/M^{3/2}+b_2/M^{5/2}+\cdots
\end{displaymath}

\begin{table}
\begin{center}
\begin{tabular}{cccc}
\hline\hline
$M$ & $K=8$ & $K=10$ & $K=12$ \\
\hline
45  &   0.036264215181387
    &   0.036264215181095
    &   0.036264215180475 \\
46  &   0.036264215181343
    &   0.036264215181088
    &   0.036264215181354 \\
47  &   0.036264215181305
    &   0.036264215181073
    &   0.036264215180915 \\
48  &   0.036264215181271
    &   0.036264215181060
    &   0.036264215180994 \\
49  &   0.036264215181240
    &   0.036264215181048
    &   0.036264215180970 \\
50  &   0.036264215181214
    &   0.036264215181038
    &   0.036264215181001 \\
51  &   0.036264215181190
    &   0.036264215181029
    &   0.036264215180972 \\
52  &   0.036264215181168
    &   0.036264215181021
    &   0.036264215180972 \\
53  &   0.036264215181149
    &   0.036264215181013
    &   0.036264215180969 \\
54  &   0.036264215181131
    &   0.036264215181006
    &   0.036264215180967 \\
55  &   0.036264215181116
    &   0.036264215181000
    &   0.036264215180962 \\
\hline\hline
\end{tabular}
\end{center}
\caption{\label{tab:SAPampl} 
Estimates of the critical SAP amplitude ${\cal P}(x_c)$.}
\end{table}
 
\noindent
Indeed, the observation that $a_1$ and $a_2$ vanish follows from the assumed
asymptotic behaviour of $p_m$.
In Table~\ref{tab:SAPampl} we have listed estimates for  ${\cal P}(x_c)$ obtained
using various values of $M$ and the cut-off $K$ in the asymptotic form. From this
we confidently estimate that ${\cal P}(x_c)=0.0362642151808(2)$.

\begin{figure}[htbp]
  \centering
  \includegraphics[scale=1.0]{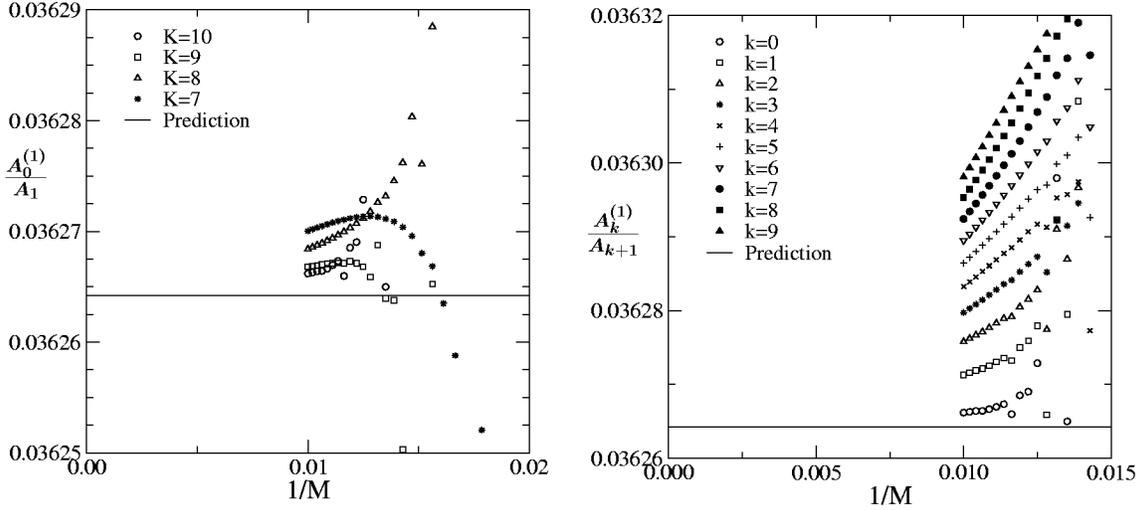}
  \caption{\label{fig:Sap1A} Estimates of the amplitude ratio $A_k^{(1)}/A_{k+1}$ for 
  self-avoiding polygons with one puncture of arbitrary size.}
\end{figure}

In figure~\ref{fig:Sap1A} we have plotted estimates of the amplitude ratio 
$A_k^{(1)}/A_{k+1}$ for self-avoiding polygons with one puncture of arbitrary size.
In the leftmost panel we look at the ratio  $A_0^{(1)}/A_{1}$ using different
cut-offs $K$. The rightmost  panel shows the ratio $A_k^{(1)}/A_{k+1}$ for 
area moments up to $k=9$ using the cut-off $K=10$. The straight line corresponds
to the expected value ${\cal P}(x_c)$, using the estimate for this quantity obtained above. 
The estimates in the leftmost panel show
some variation when plotted against $1/M$, but in the limit $M\to \infty$
the estimates appear to converge to the expected value (if the trend holds). As for minimally 
punctured SAPs we see an ever closer agreement as $K$ is increased, again indicating 
that the assumed asymptotic form is reasonable. Obviously, as seen in the rightmost panel, 
the estimates for high area moments are not as close to the expected value. However,
given the trend in these estimates we are confident in stating that our
numerical analysis is consistent with the results of Theorem~\ref{theo:arb}.
The agreement is particularly impressive bearing in mind that we analyse
series with just 42 terms.

\begin{figure}[htbp]
  \centering
  \includegraphics[scale=0.95]{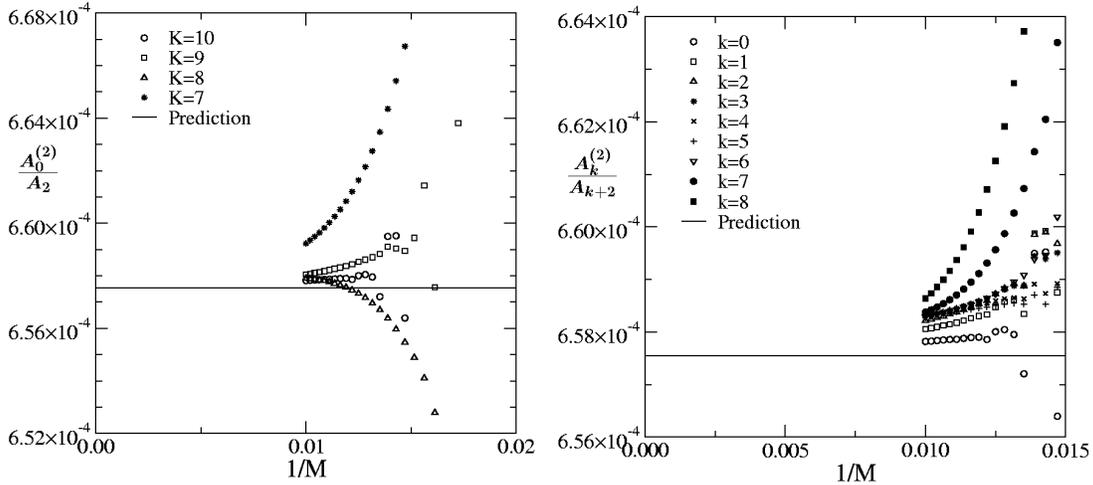}
  \caption{\label{fig:Sap2A} Estimates of the amplitude ratio $A_k^{(2)}/A_{k+2}$ 
  for self-avoiding polygons with two punctures of arbitrary size.}
\end{figure}

Finally in figure~\ref{fig:Sap2A} we have plotted our estimates of the amplitude ratio 
$A_k^{(2)}/A_{k+2}$ for self-avoiding polygons with two punctures of arbitrary size.
The straight line corresponds to the expected value ${\cal P}(x_c)^2/2$. Again
all estimates are consistent with the results of Theorem~\ref{theo:arb}, though
the numerical agreement is less convincing, but then again the series have only
38 terms.

\section{Conclusion}

We have rigorously analysed the effect of punctures on
the area law of polygon models. In particular we obtained
expressions for the leading amplitudes of the area moments 
for punctured polygons in terms of the amplitudes for unpunctured
polygons (see Theorems~\ref{theo:bounded} and~\ref{theo:arb}). For staircase
polygons this led to exact formulas for the amplitudes.
For self-avoiding polygons the formulas rely on an assumption about
the asymptotic behaviour of SAPs without punctures. They contain constants not 
known exactly but estimated numerically to a very high degree of accuracy.
Our analysis also led to conjectures
about the scaling behaviour of these models (see Conjecture~\ref{con1}).
A proof of these conjectures is an open and difficult problem. Further
numerical support of these conjectures might follow from an analysis
of critical perimeter moments along the lines of \cite{RJG04}.

The expressions for the amplitudes were thoroughly checked
numerically. For staircase polygons with up to 5 minimal punctures and
staircase polygons with one or two punctures of fixed size
we used our series expansions to find the {\em exact} generating
functions for area moments to order 10. Naturally, in all these cases the
leading amplitude agrees with the proved formulas.  Interestingly
we find that the amplitude of the correction term is proportional
to the corresponding leading amplitude with one less puncture.
For staircase polygons with one and two punctures of arbitrary size
a careful asymptotic analysis of the series for the area moments
yielded very accurate estimates for the amplitudes, again confirming
the exact formulas. Finally, we also analysed series for self-avoiding
polygons with one and two punctures (minimal as well as arbitrary).
In this case the numerical evidence is not quite as convincing, but
we did find that the numerical estimates agree with the exact formulas
to at least 3--4 significant digits.

The numerical analysis also yielded explicit expressions for correction 
terms, see Eqs.~(\ref{eq:StMB}), (\ref{eq:S1C}) and subsection \ref{sec:stair_fix}.
These correspond to contributions from puncture-boundary
interactions and from puncture-puncture interactions. It would
be interesting, but difficult to give a combinatorial proof of these results.
The difficulty of any  combinatorial proof becomes clearer when
considering the recent closed form solution of a closely related model
of punctured staircase polygons. In \cite{IJ07} one of us (IJ) considered
a model of punctured staircase polygons in which the internal polygon
is rotated by 90 degrees with respect to the outer polygon. 
The proofs in this paper never consider such restrictions on the 
placement of the internal polygon, that is, the internal polygon can 
be placed in any way one pleases. The results therefore carry over 
unaltered to the problem of rotated punctured staircase polygons. 
Interestingly, this means that the leading asymptotic form of
the coefficients is exactly the same for both models, any differences arising 
only from the correction terms. The dominant correction term to the 
half-perimeter generating function for punctured staircase polygons  
$\PGf(x)$ is $\propto  \log(1-4x)/\sqrt{1-4x}$ \cite{GJ06b}.
From the exact closed form solution to the half-perimeter generating function 
for rotated punctured staircase polygons  $\PPR(x)$ it follows that
the first correction term is  $\propto (1-4x)^{-3/4}$ \cite{IJ07}. 
These differences  indicate that combinatorial arguments for a
proof of sub-dominant behaviour must be quite subtle!

In this paper, we discussed models of punctured polygons,
where punctures are of the same type as the outer polygon.
More generally, punctures may be built from different
polygon classes. In that situation, Theorem 1 holds, with the obvious
modification, for any collection of polygon models as punctures. 
Theorem 2 holds, with the obvious modification, for any collection 
of polygon models as punctures, if the corresponding half-perimeter generating 
functions are finite at the radius of convergence $x_c$ of the half-perimeter generating
function of the outer polygon. This includes the model considered in 
\cite{IJ07} as a special case. We remark that the critical half-perimeter
generating function of the outer polygon may be infinite.

Our analytical results are based on the observation that 
puncture counting can be done using polygon area estimates, 
polygon boundary contributions being asymptotically
negligible. In particular, effects of self-avoidance 
do not influence the results. This phenomenon is also
expected to hold in higher dimension.  
Consider so-called polycubes, the generalisation of
polyominoes to three dimensions. Polycubes have been
enumerated by volume up to 18 cubes, see \cite{AB06}. 
Three-dimensional vesicles \cite{W93} are a subclass
of polycubes, having no interior holes. Let a class of three-dimensional
vesicles be given, counted by surface area, with a bounded number 
of vesicular holes. Assume that the asymptotic behaviour
of the volume moments is known for the model without
vesicular holes. If $0<\phi<1$ and the critical surface
area generating function is finite, then our method of proof 
can be adapted to describe the volume law of three-dimensional 
vesicles with vesicular holes. For the full model of closed
self-avoiding orientable surfaces of genus zero on the cubic 
lattice, however, there is numerical evidence that $\phi=1$, 
see e.g.~the review in \cite{Vanderzande}.

For models with minimal punctures, the number $p_m^{\Box(1,0)}$ 
in Eq.~(\ref{form:minprk}) also counts the number of 
polygons winding around a fixed point 
of the dual lattice. This problem and its generalisation 
to $M>1$ mutually avoiding self-avoiding polygons has been 
considered previously by Cardy \cite{C00}. It would be
interesting to consider whether this generalisation
can also be treated using the above methods. This
involves the analysis of polygon models satisfying 
$\theta=0$. If $\theta<0$, interaction terms are generally
not asymptotically negligible, so that the above
analysis does not yield asymptotically exact estimates. 

A major open question is the problem of self-avoiding polygons
with an unbounded number of punctures. Here, Theorem~\ref{theo:arb}
yields an upper bound.
Let $\QGf_{k}(x)$ denote the $k^{th}$ area moment generating function
for a model of punctured polygons with an arbitrary number of punctures.
$\QGf_{k}(x)$ is a formal power series, usually with zero radius of convergence.
The (asymptotically exact) upper bound for models with $r$ 
punctures yields an upper estimate for the number of punctured 
polygons with an arbitrary number of punctures. It is
\begin{displaymath}
\QGf_{k}(x) \ll \sum_{r=0}^\infty \frac{\PGf_{k+r}(x)}{r!}(\PGf_{0}(x))^r,
\end{displaymath}
where $\ll$ denotes coefficientwise majorisation. Let $\widetilde \QGf_{k}(x)$ 
denote the  $k^{th}$ factorial area moment generating function
for punctured polygons with an arbitrary number of punctures. 
It can be shown that an upper bound is given by
\begin{displaymath}
\widetilde \QGf_{k}(x) \ll \sum_{r=0}^\infty \frac{\widetilde \PGf_{k+r}(x)}{r!}
(\widetilde \PGf_{0}(x))^r=\left.\frac{{\rm d}^k}{{\rm d}q^k} {\cal P}(x,q)
\right|_{q={1+\cal P}(x,1)},
\end{displaymath}
where $\widetilde \PGf_{k}(x)$ is the $k^{th}$ factorial moment 
generating function, and ${\cal P}(x,q)$ is the perimeter and area 
generating function of the model without punctures.  In particular,
the half-perimeter generating function $\QGf_{0}(x)=\widetilde 
\QGf_{0}(x)$ is majorised by
\begin{displaymath}
\QGf_{0}(x) \ll {\cal P}(x,1+{\cal P}(x,1)).
\end{displaymath}

Another problem touched upon in this article is punctured
polygons in the fixed area ensemble. We gave a (non-rigorous)
argument for values of the critical
exponent in the branched polymer phase from the crossover 
behaviour of the tentative scaling function, thereby confirming 
previous results \cite{JvR92}.
In that phase, boundary effects are indeed crucial, such that 
our methods of deriving limit distributions cannot be applied in this
situation. On the other hand, area laws are expected to be
of Gaussian type, as is usually the case away from
phase transition points. The same phenomenon is expected
to occur in the fixed perimeter ensembles for $q\ne1$.

\section*{E-mail or WWW retrieval of series}

The series for the generating functions studied in this paper 
can be obtained via e-mail by sending a request to 
I.Jensen@ms.unimelb.edu.au or via the world wide web on the URL
http://www.ms.unimelb.edu.au/\~{ }iwan/ by following the instructions.

\section*{Acknowledgements}

The calculations presented in this paper would not have been possible
without a generous grant of computer time on the server cluster of the
Australian Partnership for Advanced Computing (APAC). We also used
the computational resources of the Victorian Partnership for Advanced 
Computing (VPAC). 
CR and AJG would like to acknowledge support by the German Research Council
(Deutsche Forschungsgemeinschaft) within the CRC701.
IJ and AJG gratefully acknowledge financial support from the Australian Research Council.

\section*{Appendix}

We analyse the asymptotic growth of a Cauchy product of sequences
in terms of the asymptotic growth of its constituting sequences.
The following lemma is an extension of \cite[Thm~2]{B74} to the case
of generating functions with equal radii of convergence. Its proof relies
on Lebesgue's dominated convergence theorem \cite[Thm~1.34]{R87}, which 
states conditions under which an exchange of limit and sum is allowed:  
With $n\in\mathbb N$ and $k\in\mathbb N_0$ let real numbers 
$a_{n,k}$ be given. Assume that $\lim_{n\to\infty}a_{n,k}=:a_k\in\mathbb R$ 
for all $k$ and that for all $k$ there is a bound $|a_{n,k}|\le b_k$ 
uniformly in $n\in\mathbb N$. Assume that $\sum_{k\ge0}  b_k\in
\mathbb R$. Then $\lim_{n\to\infty} \sum_{k\ge0} a_{n,k} =
\sum_{k\ge0} a_k\in\mathbb R$. 

\begin{lem}
Let two sequences $(f_n)_{n\in\mathbb N_0}$ and $(g_n)_{n\in\mathbb N_0}$
of real numbers be given, with generating functions $f(x)=
\sum_{n\ge0}f_n x^n$ and $g(x)=\sum_{n\ge0}g_n x^n$. Assume that both 
generating functions have the same positive finite radius of convergence 
$x_c$, $0<x_c<\infty$. Assume that the sequences $(f_n)$ and $(g_n)$ 
satisfy asymptotically
\begin{displaymath}
f_n \sim Ax_c^{-n}n^{\gamma-1},\qquad g_n \sim Bx_c^{-n}n^{\delta-1} \qquad
(n\to\infty),
\end{displaymath}
for nonzero numbers $A\ne0$, $B\ne0$, and for real constants $\gamma,\delta$
satisfying $\delta<0$ and $\gamma>\delta+1$. Assume that $g(x_c):=
\lim_{x\nearrow x_c} g(x)\ne0$. Then, the Cauchy product of $(f_n)$ and $(g_n)$  
satisfies
\begin{displaymath}
\sum_{k=0}^n f_{n-k} g_k = [x^n]f(x)g(x)\sim
g(x_c) f_n \qquad (n\to\infty).
\end{displaymath}
\end{lem}

\begin{proof}
Let $f_{n}:=0$ for $n<0$ and define for $n\in\mathbb N$ and for $k\in\mathbb N_0$
\begin{displaymath}
a_{n,k}:=\frac{f_{n-k}}{Ax_c^{-n}n^{\gamma-1}}g_k.
\end{displaymath}
Below, we derive a bound on $|a_{n,k}|$ uniformly in $n\in\mathbb N$ and
summable in $k\in\mathbb N_0$. Then, Lebesgue's dominated convergence theorem
can be applied to $a_{n,k}$. Since $g(x_c)\ne0$, this 
yields the statement of the lemma. 

Note first that the assumption on the asymptotic 
behaviour of $f_n$ implies the existence
of a constant $n_0\in\mathbb N$ such that for all $n\ge n_0+k$ we have the
estimate
\begin{displaymath}
\left|\frac{f_{n-k}}{Ax_c^{-n}n^{\gamma-1}}\right| 
\le 2 x_c^k \left(1-\frac{k}{n}\right)^{\gamma-1}.
\end{displaymath}
Fix such $n_0$. We distinguish the three cases $n< k$, $k\le n<n_0+k$, and $n\ge n_0+k$.
If $n< k$, we clearly have $a_{n,k}=0$. For $n\ge n_0+k$, 
the above estimate yields
\begin{eqnarray*}
\left|a_{n,k}\right|&&=\left|\frac{f_{n-k}}{Ax_c^{-n}n^{\gamma-1}}g_k\right|
\le 2\left(1-\frac{k}{n}\right)^{\gamma-1}x_c^k \left|g_k\right|.
\end{eqnarray*}
We will first consider the case $\gamma-1<0$. If $n\ge n_0+k$, we have
$1/(1-k/n)\le1+k/n_0$. This implies
\begin{eqnarray*}
\left|a_{n,k}\right| &\le 2^{2-\gamma_1} k^{1-\gamma}x_c^k\left|g_k\right|=:b_k^{(1)}.
\end{eqnarray*}
If $k\le n<n_0+k$, we estimate similarly
\begin{eqnarray*}
&\hspace{-1.5cm}\left|a_{n,k}\right|=\left|\frac{f_{n-k}}{Ax_c^{-n}n^{\gamma-1}}g_k\right|
\le |A|^{-1}\max\{1,x_c^{n_0}\} \max_{m<n_0}\{|f_m|\} (n_0+k)^{1-\gamma}x_c^k |g_k|\\
&\le |A|^{-1}\max\{1,x_c^{n_0}\}\max_{m<n_0}\{|f_m|\} \cdot
\left\{\begin{array}{cc}n_0^{1-\gamma}|g_0|, & k=0\\
(2n_0)^{1-\gamma}x_c^k k^{1-\gamma}|g_k|, & k\ne0
\end{array}
\right\}=: 
b_k^{(2)}.
\end{eqnarray*}
Define $b_k:=b_k^{(1)}+b_k^{(2)}$ for $k\in\mathbb N_0$. 
Then $|a_{n,k}|\le b_k$ uniformly in $n\in\mathbb N$, 
and $\sum_{k\ge0}b_k<\infty$.
Summability follows from
\begin{displaymath}
\sum_{k=0}^\infty x_c^k k^{1-\gamma}|g_k|<\infty,
\end{displaymath}
since by assumption $x_c^k k^{1-\gamma}|g_k|\sim |B| k^{\delta-\gamma}$ as $k\to\infty$,
where $\gamma-\delta>1$.

If $\gamma-1\ge0$, uniform estimates are obtained along the same lines.
In that situation, the factors $(1-k/n)^{\gamma-1}$ and $(n_0+k)^{1-\gamma}$
may be replaced by unity, resulting in simpler uniform bounds involving 
$\sum_{k\ge0} x_c^k |g_k|<\infty$.
\end{proof}

\end{document}